\documentclass[review]{elsarticle}

\usepackage{lineno,hyperref}
\modulolinenumbers[5]

\usepackage{tikz}

\begin{document}

\begin{frontmatter}

\title{ Numerical simulation of wetting phenomena by a meshfree particle method   }

\author[TUKL]{Sudarshan Tiwari\corref{cor1}}
\ead{tiwari@mathematik.uni-kl.de} 
\cortext[cor1]{Corresponding Author}
\author[TUKL1]{Axel Klar}
\ead{klar@mathematik.uni-kl.de}
\author[TUD]{Steffen Hardt}
\ead{hardt@csi.tu-darmstadt.de}
 
\address[TUKL]{Fachbereich Mathematik, TU Kaiserslautern,
         Gottlieb-Daimler-Strasse, 
         67663 Kaiserslautern, Germany}
\address[TUKL1]{Fachbereich Mathematik, TU Kaiserslautern,
         Gottlieb-Daimler-Strasse, 
         67663 Kaiserslautern, Germany\\
         Fraunhofer ITWM Kaiserslautern, 67663 Kaiserslautern, Germany}         
\address[TUD]{Center of Smart Interfaces, TU Darmstadt, 
 Alarich-Weiss-Str. 10, 64287, TU Darmstadt Germany}

%





\begin{abstract}
Simulations of wetting phenomena by a meshfree particle method are presented. The incompressible Navier-Stokes equations 
are used to model the two-phase flow. The continuous surface force model is used to incorporate the surface tension force.   Chorin's projection method is applied to discretize the Navier-Stokes equations. The different fluid phases are identified by assigning different colors and different material properties (density, viscosity) to the particles that remain unchanged throughout a simulation. 
Two-phase flow is captured by a one-fluid model via using 
weighted averages of the density and viscosity in a region around the fluid-fluid interface. The differential operators at each particle are computed from the surrounding cloud of particles with the 
help of the least-squares method.   
The numerical results are compared with specific analytical solutions, but also with previously considered test cases involving wetting of a container and sessile drops. A good overall agreement is found. 
 
\end{abstract}

\begin{keyword}
two-phase flow, meshfree particle method, wetting, contact angle
\end{keyword}

\end{frontmatter}

\linenumbers

\section{Introduction}
\label{section1}
 
Surface-tension driven flows occur when the surface-tension forces acting on a liquid are of equal or even larger magnitude than the inertial, viscous or gravitational forces. 
There exist numerous applications in which such flows are important, for example in the areas of microfluidics \cite{SBPH}, \cite{TCSFAHH}, coating technology \cite{WR},  
two-phase heat transfer \cite{EL} or oil recovery \cite{Babadagli}, \cite{BB}. Correspondingly, there is a high demand for efficient numerical models and schemes to compute surface-tension driven flows. 

In the past decades, a large number of CFD approaches have been presented that explicitly resolve the interface between two immiscible fluids. Often these are classified either as interface-tracking or interface-capturing schemes. In the former, the time evolution of the interface is represented by the time evolution of the numerical grid whose structure contains information about the interface shape. In the latter, the evolution of the interface is decoupled from the grid. Instead, the interface is reconstructed from field quantities represented on the grid. The most popular examples of interface capturing schemes are probably the volume-of-fluid \cite{HN}and the level-set method \cite{SSO}. While interface-tracking schemes tend to be very accurate, they are not well suited to study flows with topological changes occurring when, for example, drops or liquid sheets break up or merge. For the latter, interface capturing schemes are suitable candidates.
 
The past decades, fast progress has been made in the development of particle-based or meshfree methods to compute various types of flows. In such schemes Lagrangian particles, serving as the basic building blocks for discretization of the fluid dynamic equations, are advected with the flow. This has the advantage of a certain degree of inherent adaptivity, i.e. the numerical resolution is provided only where it is needed. Classical examples for the application of meshfree methods are astrophysical flows \cite{ML}  in which the set of Lagrangian particles co-evolves with an astrophysical structure, e.g. a plasma cloud. Often, in gas-liquid flows a similar division in an important and an unimportant  subdomain occurs. Compared to the liquid, the stresses in the gas phase are often negligible, which means that it is sufficient to compute the flow in the liquid phase alone. Meshfree methods lend themselves for this purpose. Furthermore, in contrast to interface-tracking schemes, meshfree methods easily allow studying free-surface flows with topological changes.  In total, meshfree methods appear to be ideal candidates for the simulation of complex gas-liquid flows. Among these, surface-tension driven flows form an important subclass. 

The first meshfree Lagrangian method that has been formulated  to solve fluid dynamics equations is denoted Smooothed Particle Hydrodynamics (SPH) \cite{GM}. Another meshfree Lagrangian CFD approach is the moving particle semi-implicit method \cite{KNO}. In this article we use a meshfree particle method, called Finite Pointset Method (FPM), to solve the incompressible Navier-Stokes equations. FPM is a fully Lagrangian particle method and has similar character as the SPH method except for the approximation of spatial derivatives and the treatment of boundary conditions. In SPH the spatial derivatives at an arbitrary particle position are approximated by an interpolation approach from the surrounding particles. However, in FPM the spatial derivatives are approximated using the finite-difference approach \cite{TK}, where the spatial differential operators at an arbitrary particle position are approximated by the moving least squares approach \cite{Dilts}. The Poisson equation for the pressure field is also solved in the sense of constrained least squares. The domain boundaries are represented by boundary particles, and boundary conditions are directly prescribed on those particles.   

Meshfree particle methods are  appropriate tools to simulate surface-tension driven flows. Each phase is indicated by the color of the respective particles. When particles move, they carry all the information about the flow with them such as their color, density, velocity, etc. The colors, densities and viscosity values of all particles remain constant during the time evolution. The fluid-fluid interface is easily determined with the help of the color function. To the best of the author's knowledge, the first implementation of the continuous surface force (CSF) model to account for surface-tension forces was presented in \cite{Morris}. In \cite{TK07} an implementation of the CSF model within the FPM was presented to simulate surface-tension driven flows. The present article is an extension of \cite{TK07}, devoted to studying wetting phenomena. 

The paper is organized as follows. In section \ref{model} we present the mathematical model and the numerical scheme. In section \ref{num_method} some specific aspects of the FPM are presented. The numerical test cases are presented in section \ref{num_tests}.  In some cases analytical solutions can be calculated, and the numerical solutions are compared with the analytical ones. In those cases where analytical solutions are not known, the results are compared to numerical results published earlier. 
Moreover, some convergence studies are presented in section \ref{num_tests}. The paper finishes with concluding remarks and suggestions for future work.

\section{ Mathematical model and numerical scheme}
\label{model}
\subsection{Mathematical model}
We consider two immiscible fluids, for example, liquid and gas, where both of them are incompressible. 
We use the one fluid formulation of two-phase flows from \cite{BKZ}.  We model the two-phase flows by the incompressible Navier-Stokes equations. The equations  are 
expressed in the Lagrangian form
\begin{eqnarray}
\frac{d \vec{x}}{dt} &=&\vec{v}
\label{char}
\\
\nabla \cdot \vec{ v} & = & 0
\label{conti}
\\
\frac{D \vec{ v}}{Dt} & = & - \frac{1}{\rho} \nabla \, {p}  +
\frac{1}{\rho} \nabla\cdot(2 \mu D )  + \vec{g} + \frac{1}{\rho} \vec{F}_S, 
\label{momentum}
 \end{eqnarray}
where $\vec{v}$ is the fluid 
velocity vector, $\rho$ is the density, $\mu $ is the dynamic viscosity, 
$D$ is the viscous stress tensor 
$D = \frac{1}{2} ( \nabla \vec{v} + \nabla^{T} \vec{v} ) $,   
$\vec{g}$ is the gravitational acceleration  and $\vec{F}_S$ is the surface tension force. 
In general, $\rho$ and $\mu$ are discontinuous across the interface and remain constant in each phase.  
The surface tension force $\vec{F}_S$ is 
computed using the classical continuum surface force (CSF) model \cite{BKZ}. It acts on the vicinity of the interface between the fluids. 
In the CSF model the surface tension force $\vec{F}_S$  is defined by 
\begin{equation}
\vec{F}_S = \sigma \kappa \vec{n_I}\delta_S,
\label{CSF}
\end{equation}
where $\sigma$ is the surface tension coefficient, assumed to be a constant, $\kappa$ is the curvature, 
$\vec{n}_I$ is the unit normal vector of the interface and $\delta_S$ 
is a smeared delta function, peaked at the interface.   

The equations (\ref{char} - \ref{momentum}) are solved with initial and boundary conditions. 
\subsection{Computation of the surface tension force}
Particle methods are suitable to compute the surface tension force and surface tension driven flows. The interface  
can accurately predict the flow behaviors, see \cite{Morris},\cite{TK07}. The interface can be accurately 
determined by assigning colors or flags to the particle of each phase. 
For example, we define the color $c=1$ for the gas and $c=2$ for the liquid.  The normal vector 
$\vec{n}_I$ on the interface is computed as the gradient of the color function $c$. Since $c$ is discontinuous  across 
the interface, one has to smooth it. Let $\vec{x}$ be the position of an  arbitrary particle and has neighbors with function values 
$c_j = c(\vec{x}_j)$. We smooth $c$ at $\vec{x}$ from its neighbors with the help of the Shepard interpolation rule given by 
\begin{equation}
\tilde c(\vec{x}) = \frac{\sum_{j=1}^m w_j c_j}{\sum_{j=1}^m w_j}
\label{smoothc}
\end{equation}
where $\vec{x}$ is an arbitrary particle position, $m$ is the number of neighbors and $w_j$ is the weight function given by 
\begin{eqnarray}
w_j = w( \vec{x}_j - \vec{x}; h) =
\left\{
\begin{array}{l}
\exp (- \alpha \frac{\| \vec{x}_j - \vec{x}  \|^2 }{h^2} ),
\quad \mbox{if    }  \frac{\| \vec{x}_i - \vec{x}  \|}{h} \le 1
\\  
 0,  \qquad \qquad \quad \quad \quad \quad  \mbox{else}, 
\end{array}
\right.
\label{weight}
\end{eqnarray}
where $ \alpha $ is a positive constant.  For more details we refer to section \ref{num_method}. We observe that 
the gradient of  $\tilde c$ are non-vanishing  only in a region close to the interface.  The unit normal vector is computed by 
\begin{equation}
\vec{n}_I = \frac{\nabla \tilde c}{ |\nabla \tilde c|}. 
\label{normal_I}
\end{equation}
Furthermore, the curvature is calculated using 
\begin{equation}
\kappa = -\nabla\cdot{\vec n}_I. 
\label{curvature}
\end{equation}
There exist many possible choices for $\delta_s$, but in practice, it is often approximated as 
\begin{equation}
\delta_s \approx | \nabla \tilde c|.
\label{deltaS}
\end{equation} 
We note that $\delta_s$ is non-zero in the viciinity of the interface and zero far from it. 
\subsection{Boundary conditions}
In this paper, we consider the flow in a closed container such that there are no in- and outflow boundaries. 
We have to deal only with the solid wall and interface boundary conditions. Interface boundary conditions are implicitly taken 
care by the CSF model, {\it i. e. } an explicit precription is not necessary. 
On the solid walls we use kinematic and no-slip boundary conditions, given by 
\begin{equation}
  \vec{v}\cdot \vec{n} = 0  \quad \mbox{and} \quad  \vec{v} = \vec{0},
\end{equation}
where $\vec{n}$ is the unit normal vector on a wall. 
    
In addition to the above boundary conditions, the contact angle of the liquid at the solid wall needs to be prescribed. 
The contact angle follows from the 
wetting force, which is the balance between the cohesive forces of the liquid and   the adhesive forces between 
the liquid and the wall.  The line of intersection of the three phases is called the contact line.  The 
angle between the gas liquid interface and the solid wall is called the contact angle.   In this paper we only prescribe  the static contact 
angle. The equilibrium static contact angle $\theta_s$ of a liquid drop on a  solid wall determines the wetability. 
A wetting liquid has a static contact angle less than $90^{o}$,  and a non-wetting liquid has a static contact angle larger than $90^{o}$.  
In this paper we use the method suggested by Brackbill et al \cite{BKZ}. Let $\vec{x}_w$ be a point on the solid wall with outward normal $\vec{n}$. 
Let $\vec{n}_I$ be the normal on the interface defined by (\ref{normal_I}). Before computing the curvature from (\ref{curvature}) one guarantees that the 
boundary outward normal $\vec{n}$ makes an angle $\theta_s$ with the interface normal $\vec{n}_I$. This means that the following condition must be satisfied 
\begin{equation}
\vec{n}\cdot\vec{n}_I = \cos \theta_s
\label{s_product_of_n}
\end{equation}
Therefore, we redefine the  interface normals $\vec{n}_I$ at $\vec{x}_w$ and its nearest neighbors 
within a radius $\beta ~h$ as 
\begin{equation}
\hat{\vec{n}}_I = \vec{n}\cos\theta_s + \vec{n}_{||}\sin\theta_s,
\label{corrected_n_I}
\end{equation}
where $\vec{n}_{||}$ is the unit vector parallel to the wall normal to the three-phase contact line. One can replace $\vec{n}_{||}$ by $\vec{n}_I$ computed from (\ref{normal_I}) if the real color distribution in combination with a ghost distribution is used, as described in \cite{BKZ}. The ghost distribution is obtained by reflecting the real distribution at the wall surface. The constant $\beta$ lies in the interval $(0.6,1)$ and is  problem specific. 
In our simulations we have used $\beta = 1$. 
We replace $\vec{n}_I$ by this corrected interface normal in the vicinity of  $\vec{x}_w$,  then we compute the curvature $\kappa$ from  
equation (\ref{curvature}).     
   
\subsection{Numerical scheme}
\label{num_scheme}
We consider  Chorin's projection method \cite{Chorin} in the framework of a particle method. 
Let $dt$ be a time step and set $t^n = ndt, n = 0,1,2,\ldots $. We denote, for example, ${\vec x}^n$ as the 
position of a particle at time level $n$. 
Chorin's projection scheme consists of two steps, where in the first step we compute the 
intermediate velocity $\vec{v}^{*}$ explicetely from the momentum equation without pressure term 
\begin{equation}
{\vec v}^{*} =  \vec{v}^n + \frac{dt}{\rho} \nabla\cdot(2 \mu D^{n} )  + dt \vec{g} + \frac{dt}{\rho} \vec{F}^n_S. 
\label{intermediatev}
\end{equation}
Due to the  Lagrangian formulation we do not have to deal with the nonlinear convective term. 
In the second step, called the projection step, we compute the velocity at time level $(n+1)$  
by solving the equation
\begin{equation}
\vec{v}^{n+1} = \vec{v}^{*} -  dt \; \frac{\nabla p ^{n+1}}{\rho}
\label{correctv}
\end{equation}
with the constraint that $\vec{v}^{n+1}$ satisfies the continuity equation 
\begin{equation}
\nabla\vec{v}^{n+1} = 0.
\label{constraint}
\end{equation}
In order to compute ${\vec v}^{n+1}$ we need the knowledge of $p^{n+1}$. This is 
obtained by taking the divergence of equation (\ref{correctv}) and making use of the constraint (\ref{constraint}). 
Then we get the Poisson equation for the pressure
\begin{equation}
\nabla\cdot\left(\frac{ \nabla p^{n+1}}{\rho}\right) = 
\frac{\nabla \cdot \vec{v}^{*}} {dt}.
\label{poisson}
\end{equation}
The boundary condition for $p$ is obtained by projecting equation
(\ref{correctv}) on the outward unit normal vector $\vec{n}$ at the boundary
$\Gamma$. Thus, we obtain the Neumann boundary condition
\begin{equation}
\left(\frac{\partial p}{\partial \vec{n} }\right)^{n+1} =
- \frac{\rho}{dt} (
\vec{v}^{n+1}_{\Gamma} - \vec{v}^{*}_{\Gamma}) \cdot \vec{n},
\end{equation}
where $\vec{v}_{\Gamma}$ is the value of $\vec{v}$ on $\Gamma$.
Assuming $\vec{v}\cdot\vec{n} = 0$ on $\Gamma$, we obtain
\begin{equation}
\left(\frac{\partial p}{\partial \vec{n} }\right)^{n+1} = 0 
\label{nbc}
\end{equation}
on $\Gamma$.        

In addition, we  compute the new particle positions at the $(n+1)$th level by 
\begin{equation}
\vec{x}^{n+1}  =  \vec{x}^n + dt \; \vec{v}^n. 
\end{equation}

The numerical implementation of the above scheme requires the computation of  the first and second order partial derivatives  
at every particle position.  The spatial partial derivatives at an arbitrary particle are approximated 
from its neighboring cloud of particles with the help of the weighted least squares method, described in section \ref{num_method}.  
Next, we have to solve the Poisson equation for the pressure (\ref{poisson}).
 After smoothing of the interface as described above,  the coefficient 
$1/{\rho}$ of the Poisson equation is smooth, but strongly varying near  the interface.  The same holds for the viscosity in equation (\ref{intermediatev}) for ${\vec v}^{*}$.  
 The smoothing process is similar to the smoothing of the color function.  However, 
we have to iterate the iteration process more times if the density and viscosity have higher ratios like $1000:1$ and $100:1$, 
respectively. Otherwise, the scheme becomes unstable. 
After smoothing the density, equation (\ref{poisson}) can be re-expressed as 
\begin{equation}
-\frac{\nabla\tilde\rho}{\tilde\rho}\cdot\nabla p^{n+1} + \Delta p^{n+1} = \tilde\rho\frac{\nabla\cdot{\vec v}^{*}}{dt},
\label{poisson1}
\end{equation}
where $\tilde\rho$ is the smoothed density.  
Note that, for constant density  the first term of (\ref{poisson1}) vanishes and 
we get the pressure Poisson equation. Far from the interface we have $\tilde\rho = \rho$. 
In the following section we describe the method of solving equations of type (\ref{poisson1}) by a meshfree particle method, called the Finite 
Pointset Method (FPM). 
 
\section{Finite Pointset Method (FPM)}  
\label{num_method}
FPM is a Lagrangian meshfree particle method. It has been successfully used to simulate compressible as well as incompressible flows, 
see \cite{Kuhnert}, \cite{TK} and references there in. Also it has been extended to free-surface and two-phase flows \cite{TK02}, \cite{TK07}. 
In the following subsections we present a brief description of the method.   
 
\subsection{Approximation of spatial derivatives}
\label{der_approx}
In this paper we limit ourselves to a two-dimensional spatial domain. The extension of the method to three-dimensional space is straightforward. 
Consider the computational domain $\Omega\in R^2$. Approximate $\Omega$ by particles $\vec{x}_i, i=1,\ldots,N$, whose distribution can be quite 
irregular. These particles serve as numerical grid points. 
Let $\psi(\vec{x})$ be a scalar function and $\psi_i = \psi(\vec{x}_i)$ its values for $i=1,\ldots, N$. 
We consider the problem to approximate the spatial derivatives at an arbitrary point $\vec{x} \in \{ \vec{x}_i, i = 1, \ldots, N  \}$,  in terms of the values of a set of its neighboring points.   
In order to restrict the number of neighboring points we define a  weight function $w = w(\vec{x}_i-\vec{x}, h)$ 
with small compact support of size $h$. The size of $h$ has to be chosen such that we have at least a minimum 
number of particles, for example, in $2D$, we need at least $5$ neighboring particles. In practice we define $h$ as $2.5$ to $3$ times the initial spacing of particles, keeping in mind that this is a user defined factor. 
The weight function can be quite arbitrary. In our case we consider
the Gaussian weight function defined in (\ref{weight}),  where $ \alpha $ is equal to
$6.25$. Let $ P(\vec{x}, h) = \{ \vec{x}_j :j=1,2,\ldots,m \} $ be the
set of $ m $ neighboring points of $ \vec{x} $ in a circle of radius $h$.
 
Consider $m$ Taylor expansions of $\psi (\vec{x}_i)$ around $ \vec{x} = (x,y)$
\begin{eqnarray}
\psi(x_j,y_j)= \psi(x,y)+\frac{\partial \psi}{\partial x} (x_j - x) + 
\frac{\partial \psi}{\partial y} (y_j - y) + 
\frac{1}{2} \frac{\partial ^2 \psi}{\partial x^2} (x_j - x)^2 + 
\nonumber 
\\
\quad \quad \quad \frac{\partial ^2 \psi}{\partial x\partial y} (x_j - x)(y_j-y) + 
\frac{1}{2} \frac{\partial^2 \psi}{\partial y^2} (y_j - y)^2 + e_j
\label{taylor}
\end{eqnarray}
for $j = 1, \ldots, m$,
where $ e_j $ is the residual error. 
Denote the coefficients
\begin{center}
$
a_0 = \psi(x,y), \;
a_1 = \frac{\partial\psi}{\partial x}, \;
a_2 = \frac{\partial\psi}{\partial y}, \; $\\
$
a_3 = \frac{\partial^2\psi}{\partial x^2}, \;
a_4 = \frac{\partial^2\psi}{\partial x\partial y}, 
a_5 = \frac{\partial^2\psi}{\partial y^2}. \;
$
\end{center}
Note that $a_0$ is known, so we have five unknowns $a_i, i = 1,\ldots,5$.  
Now we have to solve $m$ equations for five unknowns . For $m > 5$ this system is 
overdetermined and can be written in matrix form as 
\begin{equation}
\vec{e}= M \vec{a} -  \vec{b},
\label{error}
\end{equation}
where $M=$
\begin{eqnarray}
\left( \begin{array}{ccccc}
dx_1 & ~dy_1 & ~\frac{1}{2}dx^2_1 & ~dx_1 dy_1 & ~\frac{1}{2} dy^2_1    \\
\vdots  &\vdots & \vdots &\vdots &\vdots   \\
dx_m & ~dy_m & ~\frac{1}{2}dx^2_m  & ~dx_m dy_m & ~\frac{1}{2} dy^2_m 
 \end{array} \right)
\label{matrixM}
\end{eqnarray}
$ \vec{ a} = \left ( a_1 , a_2 , \ldots   a_{5} \right )^T , \;
\vec{ b} =  \left ( \psi_1 - a_0, \ldots , \psi_m - a_0  \right )^T $,
$\vec{ e} = \left ( e_1, \ldots , e_m \right )^T $ 
and
$dx_j = x_{j} - x, \;  dy_j = y_{j}-y$. 

The unknowns $a_i$ are computed by minimizing a weighted error over
the neighboring points. Thus, we have to minimize the following quadratic form
\begin{equation}
J = \sum_{i=1}^{m} w_i e_i^2  = (M \vec{a} - \vec{b})^T W (M \vec{a} -\vec{b}), 
\label{functional}
\end{equation}
where 
\begin{eqnarray*}
W=\left( \begin{array}{cccc}
w_1 & 0 & \cdots& 0 \\
\vdots & \vdots & \cdots & \vdots \\
0 & 0 & \cdots & w_m  
\end{array} \right).
\end{eqnarray*}
The minimization of $ J $ with
respect to $\vec{a}$ formally yields ( if $M^T W M$ is nonsingular)
\begin{equation}
\vec{a} = (M^T W M)^{-1} (M^T W) \vec{b}.
\label{lssol}
\end{equation}

\subsection{Particle method for solving the Poisson equation}
 
With equation (\ref{poisson1})  we have to solve a linear 
partial differential of second order of the form 
\begin{equation}
\vec{B}\cdot\nabla\psi + C\Delta\psi = f,
\label{ellipeq}
\end{equation}
where $\vec{B},C $ and $f$ are given. 
The equation is solved 
with Dirichlet or Neumann boundary conditions 
\begin{equation}
\psi = g \quad \quad \quad \mbox{or} \quad 
\frac{\partial\psi}{\partial\vec{n}} = \phi. \label{ellipeqbc}
\end{equation} 

In fact, we can substitute the partial differential operators appearing in equation (\ref{ellipeq}) by the components of $a$ from equation (\ref{lssol}). This approach was first proposed in \cite{LO}. However, 
it has some difficulties when dealing with the Neumann boundary condition. 
In the following we describe a meshfree particle method, initially proposed 
in  \cite{TK1}, which is more stable compared to the method given in \cite{LO} (see \cite{IT} for details). 
Moreover, the method can easily handle the Neumann boundary condition and has a second-order convergence. 
 
We consider again an arbitrary particle  position $(x,y)$ having $m$ neighbors, as in subsection \ref{der_approx}. 
We reconsider the $m$ Taylor expansions of equation (\ref{taylor}).  We add the constraint that at particle position $(x,y)$ the partial differential equation (\ref{ellipeq}) 
should be satisfied. If the point $(x,y)$ lies on the boundary, also the boundary condition (\ref{ellipeqbc}) needs to be satiesfied. Therefore,  we add the equations 
(\ref{ellipeq}) and (\ref{ellipeqbc}) to these $m$ equations (\ref{taylor}). Equations (\ref{ellipeq}) and (\ref{ellipeqbc}) are re-expressed as 
\begin{eqnarray}
\label{constraint1}
B_1 a_1 + B_2 a_2 + C (a_3 + a_5 ) = f\\
\label{constraint2}
n_x a_1 + n_y a_2  = \phi,
\end{eqnarray}
where $n_x, n_y$ are the $x,y$ components of the unit
normal vector $\vec{n}$ on the boundary $\Gamma$.
 
In this formulation the function values $\psi(x,y)$ are not known a priori, therefore,  the term $a_0$ is also unknown. 
In total we have six unknowns $a_i, i = 0,1,\ldots, 5$. For the interior particles we add equation (\ref{constraint1}), and 
for boundary particles with Neumann boundary conditions equation (\ref{constraint2}) as constraints.  
Now we have to solve $m+1$ equations for six unknowns.  For $m+1 > 6$ this system is 
overdetermined with respect to the unknowns $a_i$ and can be written 
in matrix form (\ref{error}),  where the matrix $M$ differs from (\ref{matrixM}) and is given by 
\begin{eqnarray}
\left( \begin{array}{cccccc}
 1 & ~dx_1 & ~dy_1 & ~\frac{1}{2}dx^2_1 & ~dx_1 dy_1 & ~\frac{1}{2} dy^2_1    \\
\vdots & \vdots  &\vdots & \vdots &\vdots &\vdots   \\
1  &~dx_m & ~dy_m & ~\frac{1}{2}dx^2_m  & ~dx_m dy_m & ~\frac{1}{2} dy^2_m \\
0 & ~B_1 & ~B_2 & ~C & ~0 & ~C \\
0 & ~n_x  & ~n_y & ~0  &~0  &~0
 \end{array} \right), 
 \label{matrixM1}
\end{eqnarray}
where $\vec{ a} = \left ( a_0, a_1  , \ldots   a_{5} \right )^T , \;
\vec{ b} =  \left ( \psi_1 , \ldots , \psi_m, f, g \right )^T $ and
$\vec{ e} = \left ( e_1, \ldots, e_m, e_{m+1}, e_{m+2} \right )^T $. 
From a programming point of view, we set $n_x = n_y = 0$ for the interior particles. 
For the Dirichlet boundary particles, we directly prescribe the boundary conditions, and for the 
Neumann boundary particles we set  $B_1 = B_2 = C = 0$ and $f=0$. 

Similarly, the unknowns $a_i$ are computed by minimizing a weighted error function and obtained in the form (\ref{lssol}).  
In (\ref{lssol}) the vector $( M ^T  W)\vec{ b}$ is explicitely given by
\\
\begin{eqnarray}
( M ^T  W)\vec{ b} =
\left( \sum_{j=1}^m w_j \psi_j , \;
\sum_{j=1}^m w_j dx_j \psi_j + B_1 f + n_x \phi,
\right.  
\nonumber
\\ \left.
\sum_{j=1}^m w_j dy_j \psi_j + B_2 f + n_y \phi, \;  
\frac{1}{2}\sum_{j=1}^m w_j dx^2_j \psi_j + C f , \;
\right.
\nonumber 
\\ \left.
\sum_{j=1}^m w_j  dx_j dy_j   \psi_j,  \;
\frac{1}{2}\sum_{j=1}^m w_j dy^2_j \psi_j + C f \;
\right )^T.
\end{eqnarray}

Equating the first components on both sides of equation (\ref{lssol}), we get 
\begin{eqnarray}
\psi = Q_{1} \left(\sum_{j=1}^m w_j \psi_j  \right) +
Q_{2} \left ( \sum_{j=1}^m w_j dx_j \psi_j + B_1 f + n_x \phi \right) +
\nonumber
\\
Q_{3}\left(\sum_{j=1}^m w_j dy_j \psi_j + B_2 f + n_y \phi \right) +
Q_{4} \left( \frac{1}{2}\sum_{j=1}^m w_j dx^2_j \psi_j + C f\right) +
\nonumber 
\\
Q_{5} \left(  \sum_{j=1}^m w_j dx_j dy_j \psi_j \right) +
Q_{6} \left( \frac{1}{2}\sum_{j=1}^m w_j dy^2_j \psi_j + C f\right), 
\end{eqnarray}
where $Q_{1}, Q_{2}, \ldots, Q_{6}$ are the components of the first row of
the matrix $( M^T  W  M)^{-1}$.
Rearranging the terms, we have
\begin{eqnarray}
\nonumber
\psi - \sum_{j=1}^m w_j\left ( Q_{1} + Q_{2} dx_j +
Q_{3} dy_j +  Q_{4} \frac{dx^2_j}{2} +
Q_{5} dx_j ~dy_j + Q_{6} \frac{dy^2_j}{2}   \right) \psi_j =
\\ 
\left ( Q_{2} B_1 + Q_{3} B_2  + Q_{4} C
+ Q_{6} C  \right ) f + 
\left(Q_{2} n_x + Q_{3} n_y  \right ) \phi. \quad \quad \quad
\label{sparsesy}
\end{eqnarray}

Equation (\ref{sparsesy}) is  for an arbitrary particle $\vec{x}$, which is one of the particle $\vec{x}_i, i = 1, \ldots , N$, having $m(i)$ neighbors at    
$\vec{x}_{i_j}$.  We repeat the computation of equation (\ref{sparsesy}) for all 
particles $i=1,\ldots,N$, giving the following sparse linear 
system of equations for the unknowns $\psi_i, i=1,\ldots, N$
\begin{eqnarray}
\nonumber
\psi_i - \sum_{j=1}^{m(i)} w_{i_j}\left ( Q_{1} + Q_{2}  dx_{i_j} +
Q_{3} dy_{i_j}  + Q_{4} \frac{dx^2_{i_j}}{2} +
Q_{5} dx_{i_j} dy_{i_j} +
Q_{6} \frac{dy^2_{i_j}}{2} \right ) \psi_{i_j} = 
\\ 
\left ( Q_{2} B_1 + Q_{3} B_2  + Q_{4} C
+ Q_{6} C  \right ) f_i + 
\left(Q_{2} n_x + Q_{3} n_y  \right ) \phi_i.  \quad \quad \quad
\label{sparsesy1}
\end{eqnarray}
In the matrix form we have
\begin{equation}
\label{sparsesy2}
L\vec{\Psi} = \vec{R},
\end{equation}
where $\vec{R}$ is the right-hand side vector, $\vec{\Psi}$ is the unknown vector and 
$L$ is the sparse matrix having non-zero entries only for neighboring particles.  

The sparse system (\ref{sparsesy2}) can be solved by some iterative
methods. In this paper we apply the method of Gauss-Seidel.  
In the projection scheme it is also necessary to prescribe  initial values for the pressure 
at time $t = 0$. For example, we can prescribe a vanishing pressure initially. Then, 
in the time iteration the initial values of the 
pressure for time step 
$n+1$ are taken as the values from time step $n$.
Usually, solving the pressure Poisson equation will require more iterations in the first few time steps. After a certain number of time steps, the pressure values at the old 
time step are close to those of new time step, so the number of iterations required gets reduced.
 
The iteration process is stopped if the relative error satisfies
\begin{equation}
\frac{\sum_{i=1}^N |\psi_i^{\tau + 1} - \psi_i ^{(\tau )} | }
{\sum_{i=1}^N |\psi^{(\tau + 1)}_i |} < \epsilon, 
\label{error}
\end{equation}
where $\tau = 0, 1, 2, \ldots $, and the approximation to the solution is defined by 
$\psi(\vec{x}_i): = \psi^{(\tau +1)}(\vec{x}_i), i = 1, \ldots, N $. 
The parameter $ \epsilon $ is a small positive constant and can be 
defined by the user. The required number of iterations depends on the values of 
$\epsilon$ and $h$.   


 \section{Numerical tests}
\label{num_tests}

\subsection{Elliptic equation}
As a first numerical test we study the diffusion equation 
\begin{equation}
\nabla \cdot \left(k\nabla \psi \right ) = f  \quad \mbox{in} \quad  \Omega, 
\label{diff_eqn}
\end{equation}
where $k$ is  smooth, but strongly varying near the interface.  
Away from the interface we have 
\begin{eqnarray}
k =
\left\{
\begin{array}{l}
 k_1 \; \mbox{in}  \;\Omega_1
\\ \nonumber
 k_2 \; \mbox {in} \; \Omega_2 
\end{array}
\right.
\end{eqnarray}
with $k_1\ne k_2$ and $\Omega = \Omega_1 \cup  \Omega_2$.  
Equation (\ref{diff_eqn}) is similar to the pressure Poisson equation in the projection step, therefore, the accuracy  
of its solution is very important for our numerical scheme. For the sake of simplicity, we consider Dirichlet boundary 
conditions. In the following we study the numerical solution for 
two examples with different interfaces between the two subdomains, compare \cite{IT}. We  monitor the  error between the numerical solution of the above problem with smoothed coefficient and the  exact solution of the problem with discontinuous coefficients given by the above two  constants in the respective domains. For the latter  problem we have to prescribe additional conditions at the interface, i.e. 
 the continuity of the solution and normal component of the flux through the interface.

\subsubsection{Example 1}
\begin{figure}
\centering
\includegraphics[keepaspectratio=true, width=.5\textwidth]{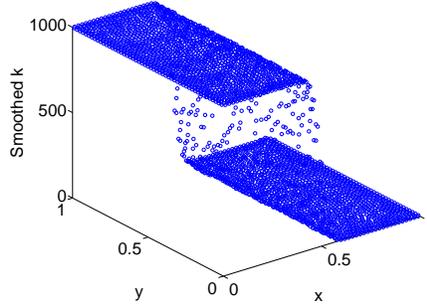}
\caption{Smoothed value of diffusion coefficient $k$ }
\label{smoothed_k}
\end{figure}

\begin{figure}
\includegraphics[keepaspectratio=true, width=.45\textwidth]{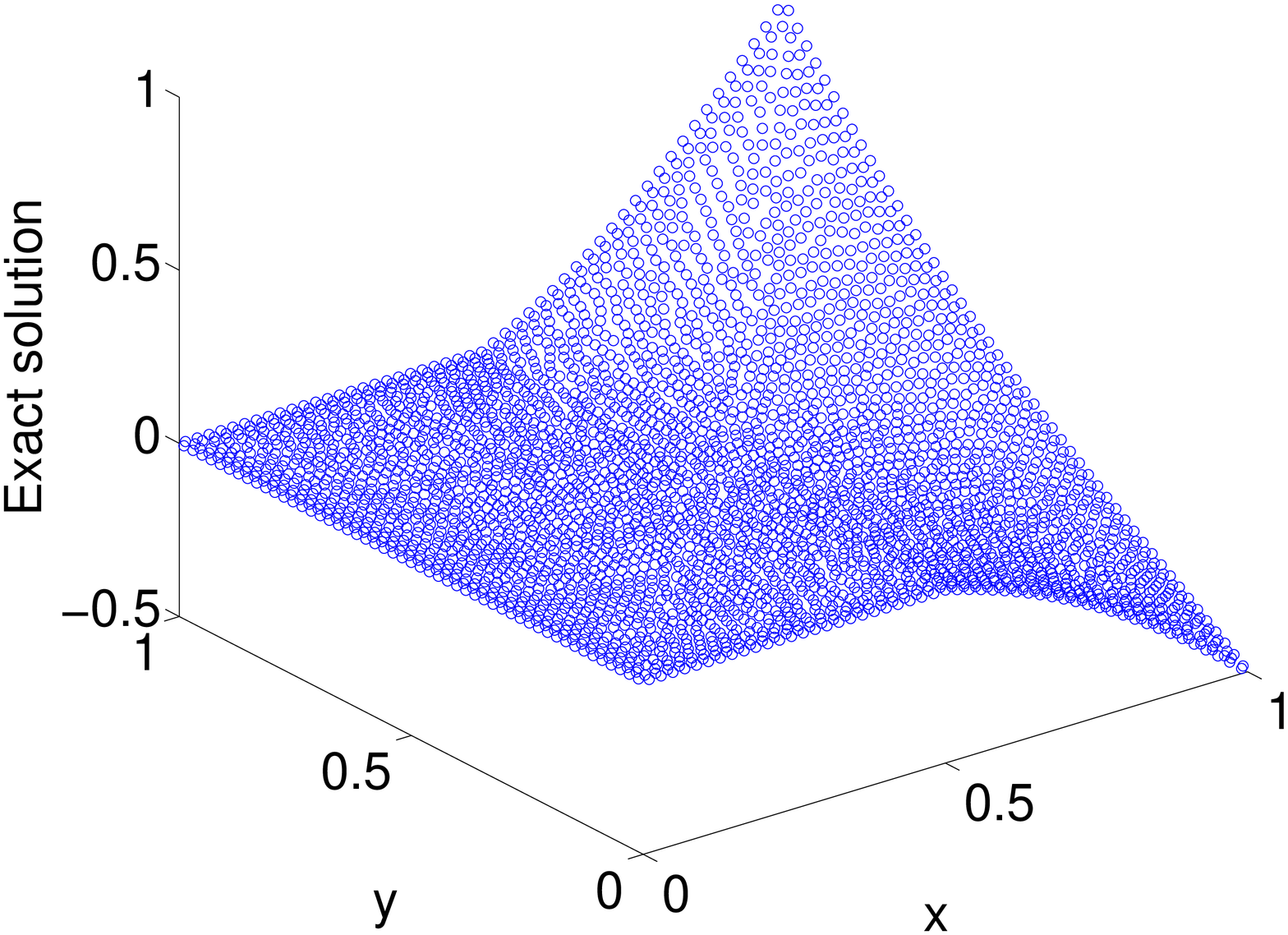}
\hfill
\includegraphics[keepaspectratio=true, width=.45\textwidth]{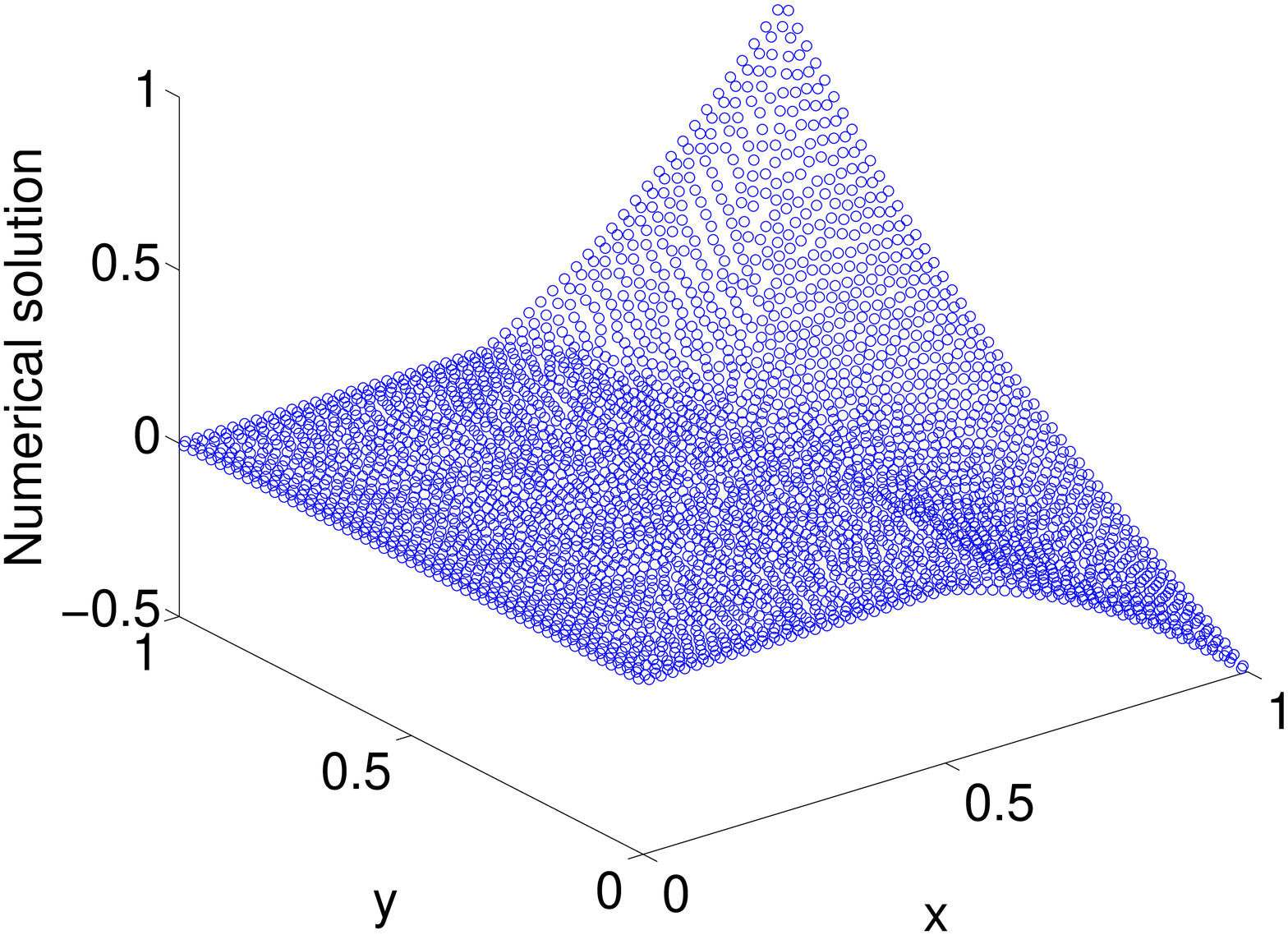}
\caption{Left: exact solution and right: numerical solution for example 1.}
\label{compare_elliptic}
\end{figure}

Let $\Omega = [0, 1]\times [0, 1]$ an unit square and decompose it into two domains 
with  interface $x = 0.5$. 
Consider $k = k_1 = 1000$ for  $x \le 0.5$ and $k = k_2  = 1$ for $x > 0.5 $. 
Using  $k = k_1$ and $k=k_2$ respectively, we define in the two domains
\begin{equation}
\varphi^{}(x,y) = \frac{1}{k} \sin (\frac{\pi x}{2})(x-\frac{1}{2})(y-\frac{1}{2})(1+x^2 + y^2). 
\label{exact1}
\end{equation}
Then, we consider a source term $f$ given in each of the domains by  
\begin{eqnarray}
f = \nabla \cdot \left(k\nabla \varphi \right ).
 \end{eqnarray}
The Dirichlet boundary conditions are   given by the values of $\varphi$ on the boundary. Obviously, the exact solution of  (\ref{diff_eqn}) with discontinuous  coefficients given by the two  constants in the respective domains  is given by 
\begin{equation}
\psi = \varphi. 
\label{exact1}
\end{equation}
We note that   solution and  flux are continuous at the interface. 
This solution is compared to the numerical solution of  (\ref{diff_eqn})
with a smoothed diffusion coefficient. 
In Fig. \ref{smoothed_k} we have plotted the 
smoothed values of  $k$.
 We iterate the smoothing process three times 
using the Shepard interpolation rule (\ref{smoothc}). 
In Fig. \ref{compare_elliptic} we show the exact and numerical solutions for $h=0.04 $, which corresponds to a 
total number of particles equal to $3417$.  
Furthermore, we have performed a convergence study, see Table \ref{table1}, where we plot the maximum error between 
the exact solution with discontinuous  $k$ and the numerical solutions with smoothed $k$. We observe that  the order of convergence is approximately one. 
In this case, using  the smoothing of the coefficent near the interface reduces   the order of convergence from two to one.  

\begin{table}
\begin{center}
\begin{tabular}{|r|r|r|}
\hline
$h$ & $N$ & $L^{\infty}$ error   \\ \hline 
$0.08$ & $909$      & $1.284\times 10^{-1}$  \\  
$0.04$ & $3222$    & $6.7518\times 10^{-2}$   \\
$0.02$ &$13177$   & $3.5763\times 10^{-2}$   \\ 
$0.01$  &$52089$  & $1.8636\times 10^{-2}$   \\ \hline
\end{tabular}
\caption{Convergence study for example 1}
\label{table1}
\end{center}
\end{table}

\subsubsection{Example 2}
Again we consider the unit square as a computational domain. We again decompose the domain into two parts. However,  
the interface is now defined as an ellipse
\begin{equation}
 (x-\frac{1}{2})^2 + 4(y-\frac{1}{2})^2 = r^2,
\end{equation}
where $r = 0.1$. We consider $k=1000$ inside the ellipse and $k=1$ elsewhere. 
The manufactured solution for the problem with discontinuous diffusion coefficient is  given by
\begin{equation}
\psi(x,y) = \varphi(x,y) = \frac{1}{k}\sin(\frac{\pi x}{2}) \left[ (x-\frac{1}{2})^2 + 4(y-\frac{1}{2})^2 - r^2 \right ] (1 + x^2 + y^2). 
\end{equation} 
 The maximum error between the 
exact solution with discontinuous coefficient and the numerical solutions of the problem with smoothed coefficient can be found in Table \ref{table2}.  The errors in Table \ref{table2} show again  first order convergence. 

\begin{table}
\begin{center}
\begin{tabular}{|r|r|r|}
 \hline
$h$ & $N$ & $L^{\infty}$ error  \\ \hline 
$0.08$ & $909$      &  $3.272\times 10^{-1} $ \\ \hline
$0.04$ & $3222$    & $1.186 \times 10^{-1} $  \\ \hline
$0.02$ &$13177$   & $4.3195\times 10^{-2}$ \\ \hline
$0.01$  &$52089$  & $1.7366\times 10^{-2}$  \\ \hline
\end{tabular}
\caption{Convergence study for example 2}
\label{table2}
\end{center}
\end{table}

\subsection{Flow induced by wall adhesion} 
 
The following test cases of flows induced by wall adhesion have been analyzed by Brackbill et al {\cite{BKZ} and later by Liu et al \cite{LKO}. 
We consider a shallow pool of water located at the bottom and  gas at the top of a two-dimensional tank of size $[0,0.112] \times[0,0.152]$. 
Two different cases for the equilibrium contact angle are studied, $\theta_s = 175^{o}$ and $5^{o}$. $\theta_s = 5^{o}$ corresponds to the wetting situation and 
$\theta_s = 175^{o}$ to a non-wetting one.  The following material parameters are chosen for both cases: $\rho_l = 1000 ~\mbox{kg m}^{-3}$, $\mu_l = 0.0091$ Pa s for the liquid, 
$\rho_g = 1 ~\mbox{kg m}^{-3}$, $\mu_g = 1.86\times 10^{-5}$ Pa s for the gas, and $\sigma = 0.072 ~\mbox{Nm}^{-1}$. The external forces, for example, the gravity force, are set to zero. The initial spacing of the 
particles is approximately $\Delta x \approx h/3$, where $h=0.004 $, which gives an initial total number of particles equal to $5737$. 
A fixed time step $dt = 0.0002 $ is chosen.  During a time step, we have to add particles if they leave a void and remove them if they are 
very close to each other.  If two particles are very close to each other, we replace the new one at the mean position and delete the two closed particles. The fluid 
quantities are assigned to newly added particles based on a least squares approximation from neighboring particles. 
We refer to \cite{TK02} for the algorithm of adding and removing particles.  The total number of particles remains approximately the same throughout a simulation. This so-called particle management is needed for all time dependent flow problems. No-slip boundary conditions are applied at all solid walls.
The initial contact angle of the liquid with the solid wall is $90^{o}$. If the prescribed contact angle $\theta_s$ is different from the initial contact angle, the contact line moves and the 
liquid surface deforms to achieve the contact angle $\theta_s$. 

For $\theta_s = 175^{o}$ we initialize the region below the line $y=0.02$ as liquid, the rest as gas, see Fig. \ref{pool_theta_175}(a), where blue (or dark gray) particles 
 are liquid particles and red particles (or light gray) are gas particles. In Fig. \ref{pool_theta_175} we have plotted the time evolution of the phase distribution. 
We observe that at around $t=1.9s$ the liquid detaches from the solid wall after which it moves inside the gas. Then it forms an oscillating drop. 
The snapshots are comparable with the results published in the papers by  Brackbill et al \cite{BKZ} and Liu et al \cite{LKO}.

   
 \begin{figure}
\includegraphics[keepaspectratio=true, width=.5\textwidth]{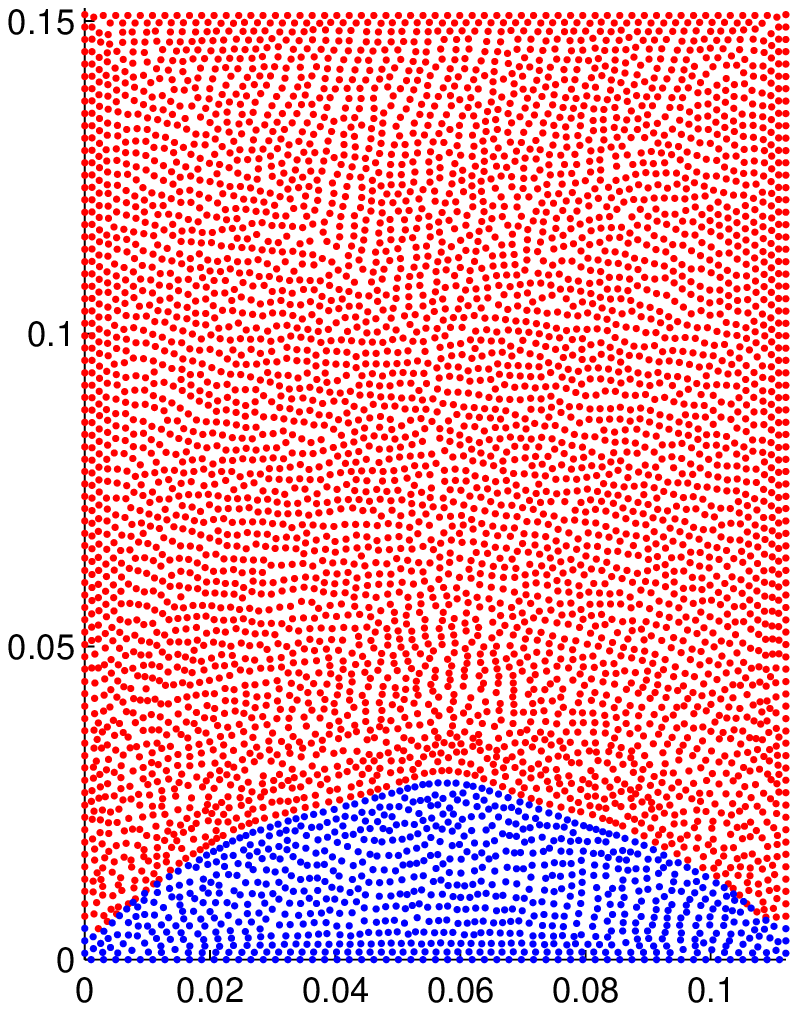}
\includegraphics[keepaspectratio=true, width=.5\textwidth]{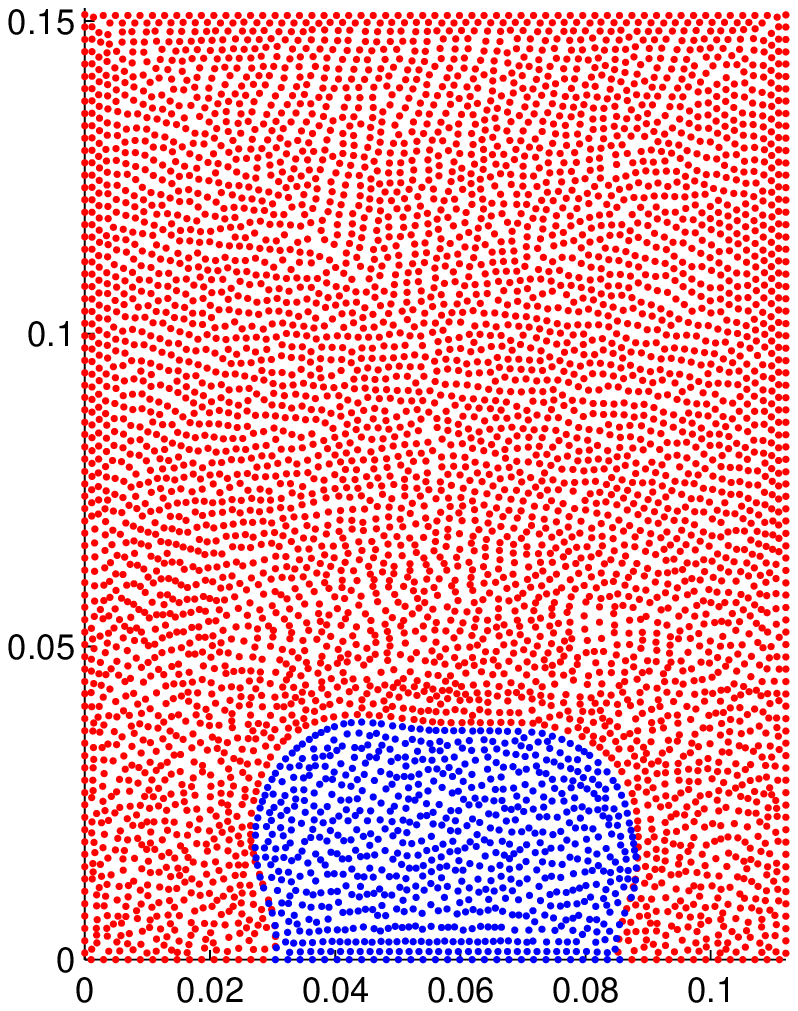}
\\
\includegraphics[keepaspectratio=true, width=.5\textwidth]{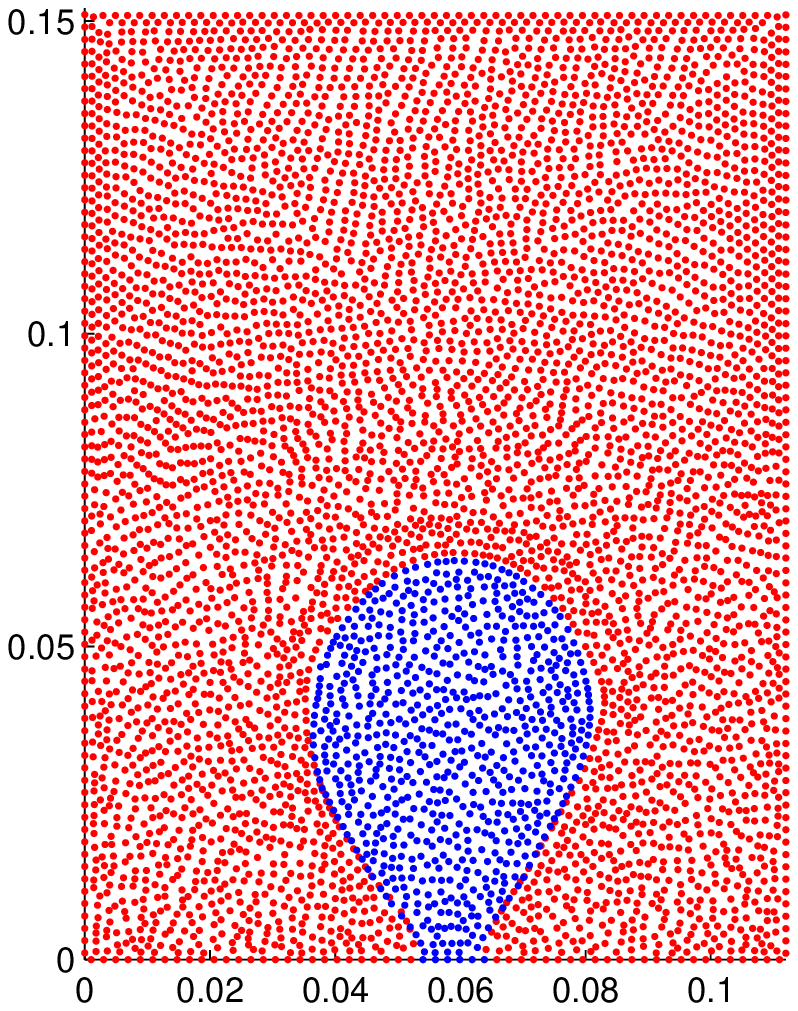} 
\includegraphics[keepaspectratio=true, width=.5\textwidth]{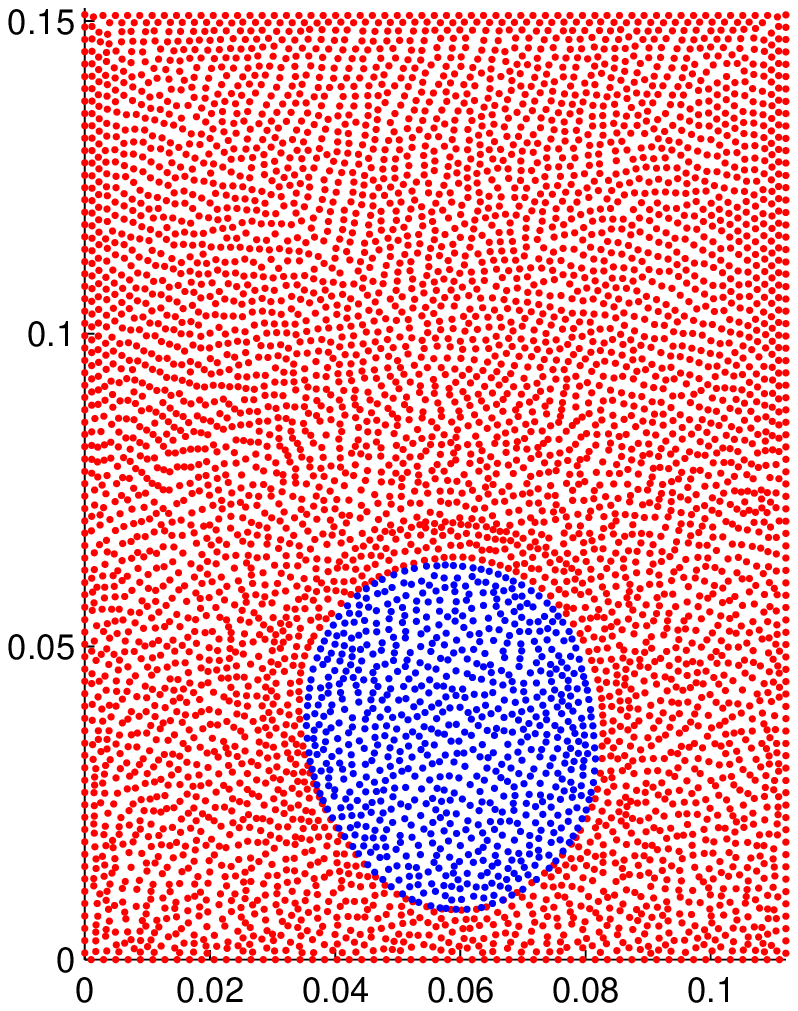}
\\
\includegraphics[keepaspectratio=true, width=.5\textwidth]{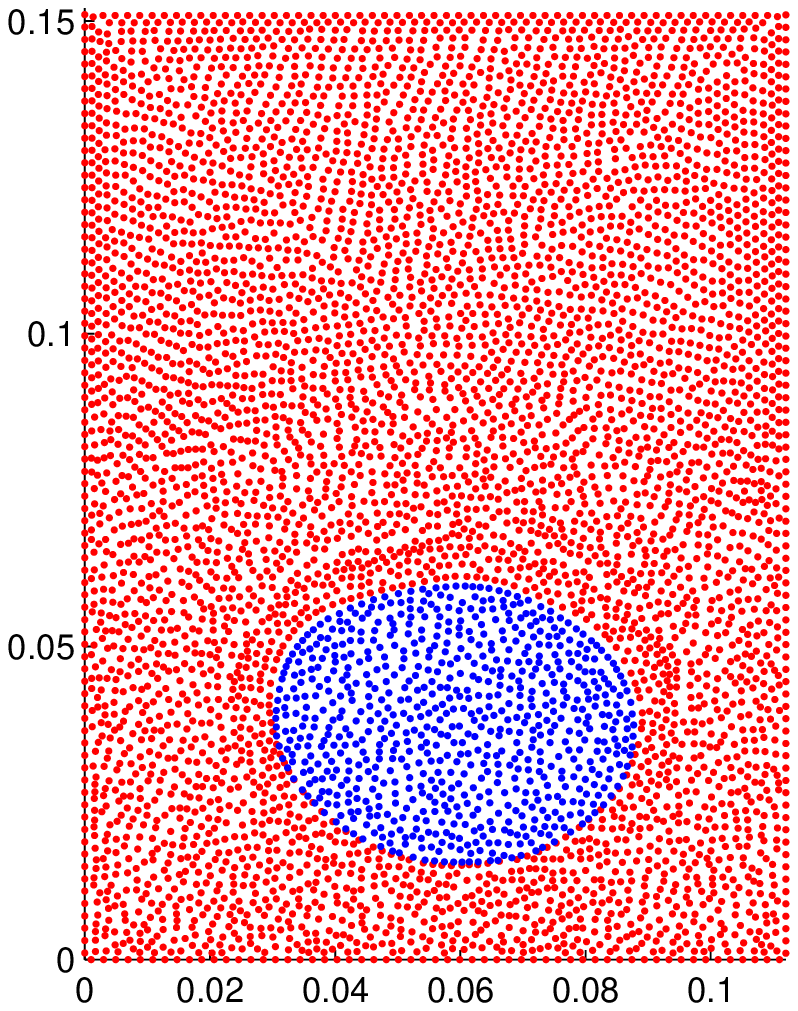} 
\includegraphics[keepaspectratio=true, width=.5\textwidth]{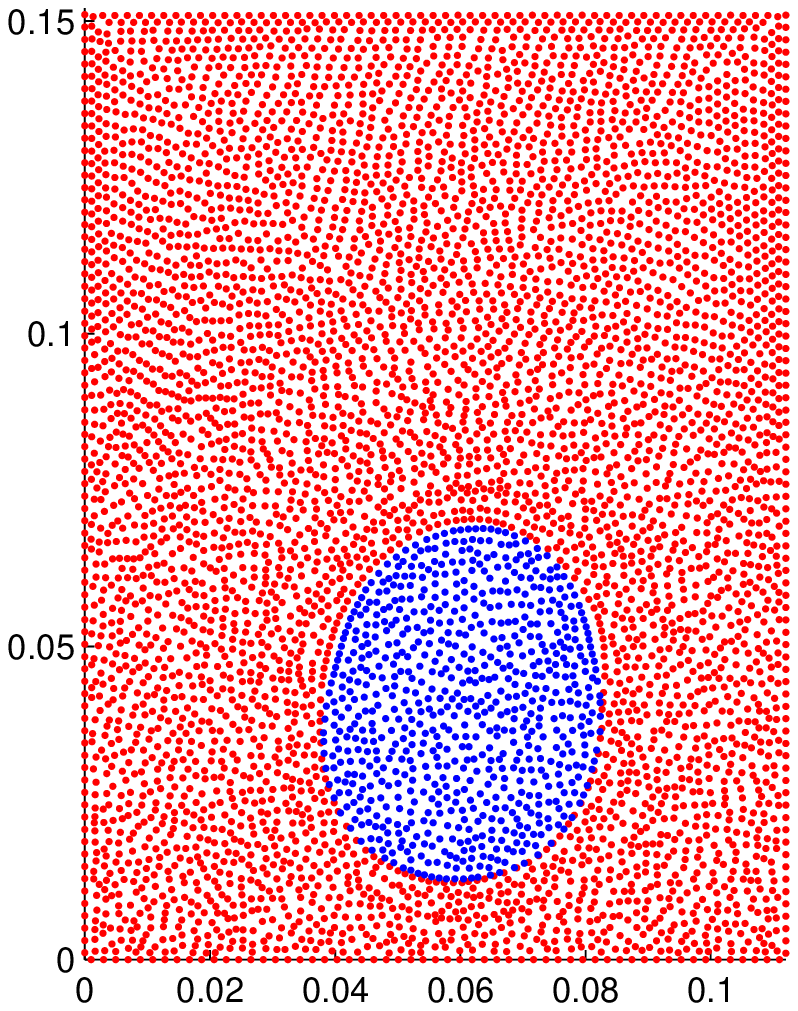}
 \caption{ Time evolution of the two phases for $\theta_s = 175^{o}$.  From left to right, first row: $t=0.5$ s and $1.0$ s, second row: $t=1.8$s and $2.0$ s, 
 third row: $t=2.4$s and $3.0$ s.  Blue (or dark gray)  particles indicate liquid particles and red (or light gray)  ones indicate gas particles. }
  \label{pool_theta_175}
\end{figure}

Next, we consider the case of $\theta_s = 5^{o}$. In this case we choose the region below the line $ y = 0.05$ to be the liquid domain, the rest of the domain is filled with gas. 
The time evolution of the two phases is plotted in Fig. \ref{pool_theta_5}. We observe that the liquid starts wetting the upper parts of the two side walls. 
At a time around $t=2.4s$ the liquid reaches its maximum height and then starts oscillating. A similar behavior was observed in \cite{BKZ}. 

 \begin{figure}
 \includegraphics[keepaspectratio=true, width=.5\textwidth]{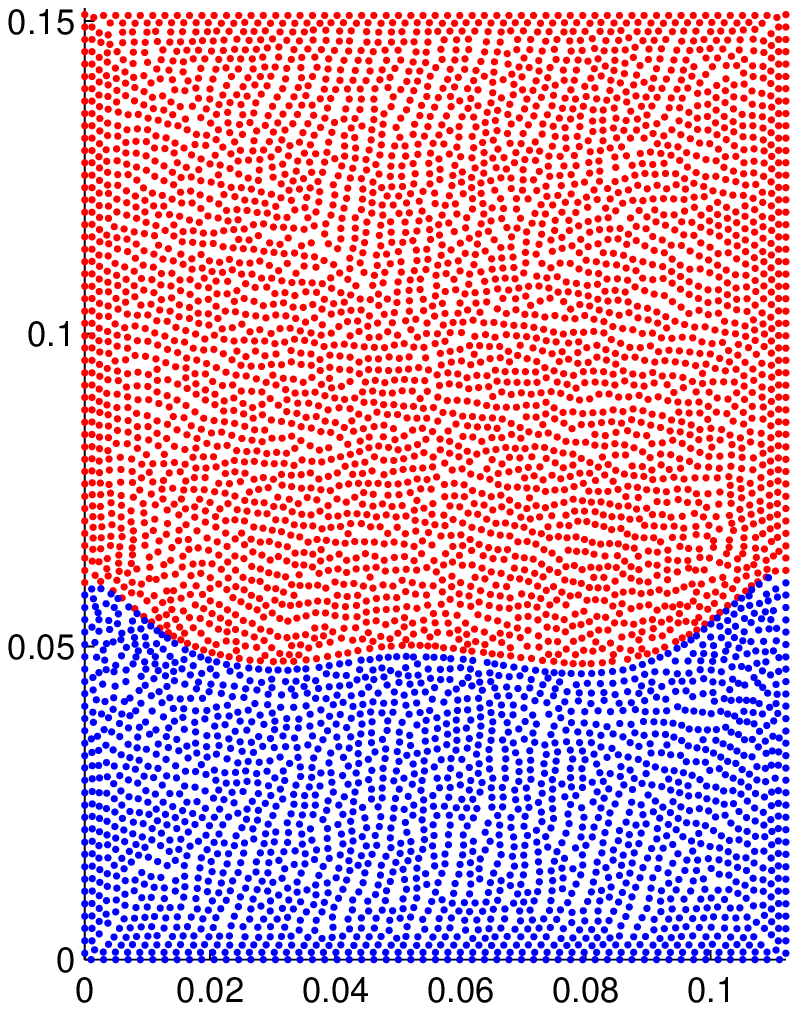}
 \includegraphics[keepaspectratio=true, width=.5\textwidth]{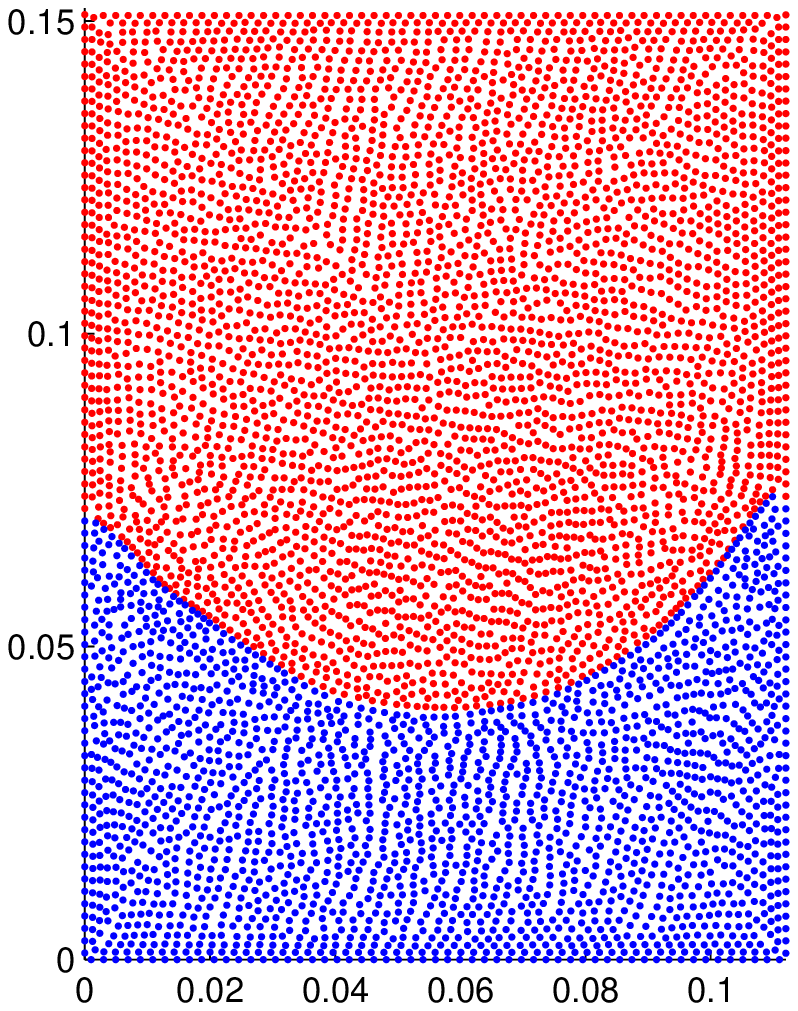}
 \\
 \includegraphics[keepaspectratio=true, width=.5\textwidth]{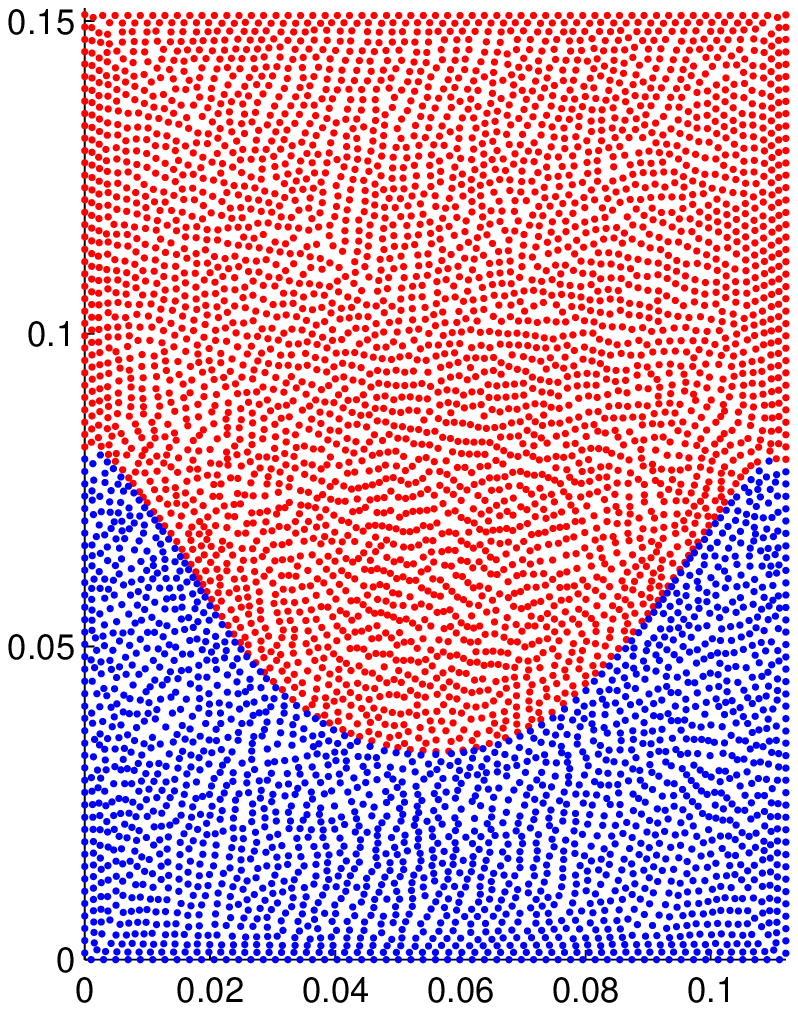}
\includegraphics[keepaspectratio=true, width=.5\textwidth]{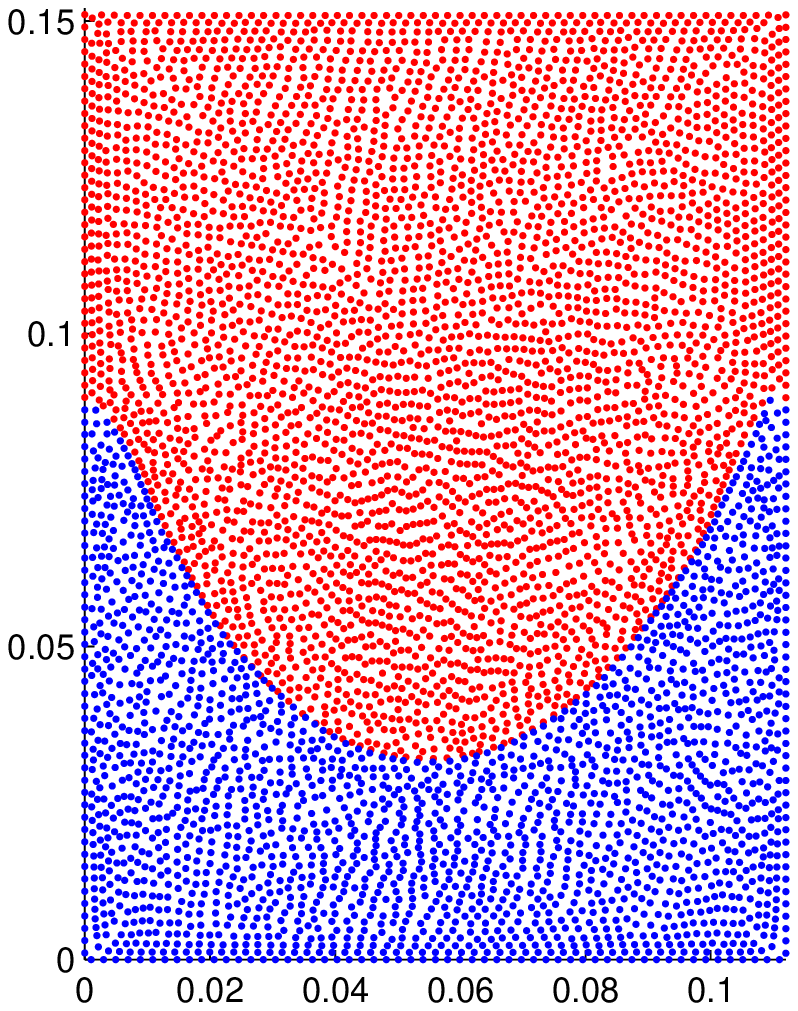}
\\
\includegraphics[keepaspectratio=true, width=.5\textwidth]{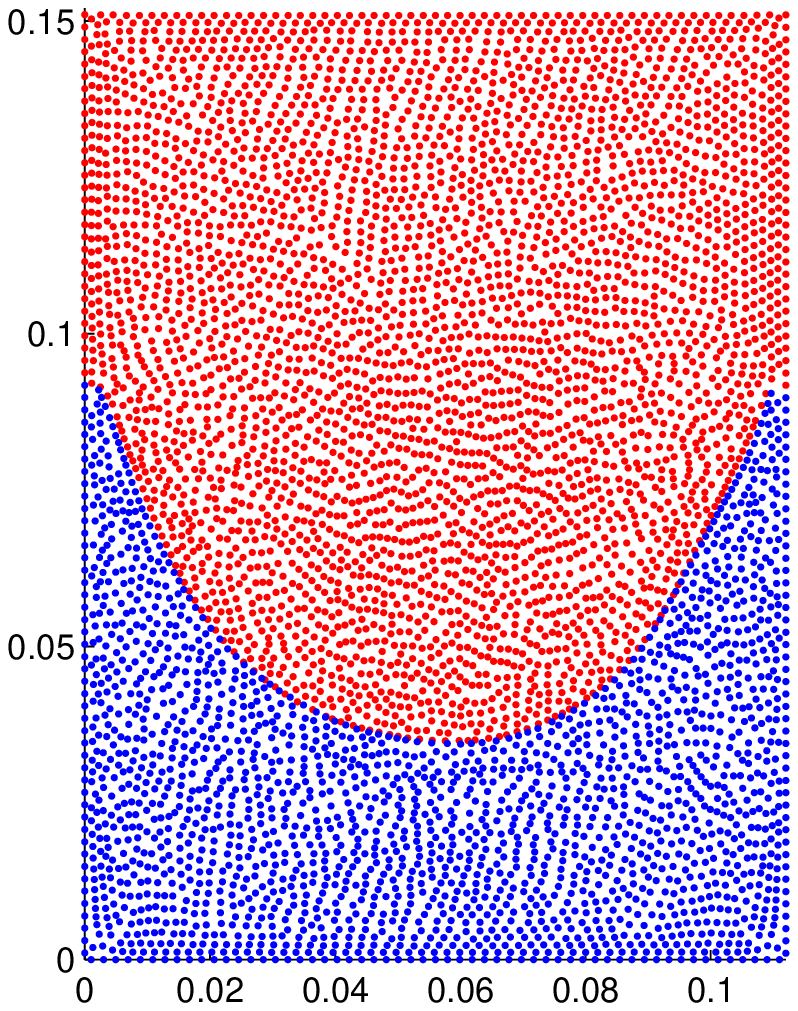}
\includegraphics[keepaspectratio=true, width=.5\textwidth]{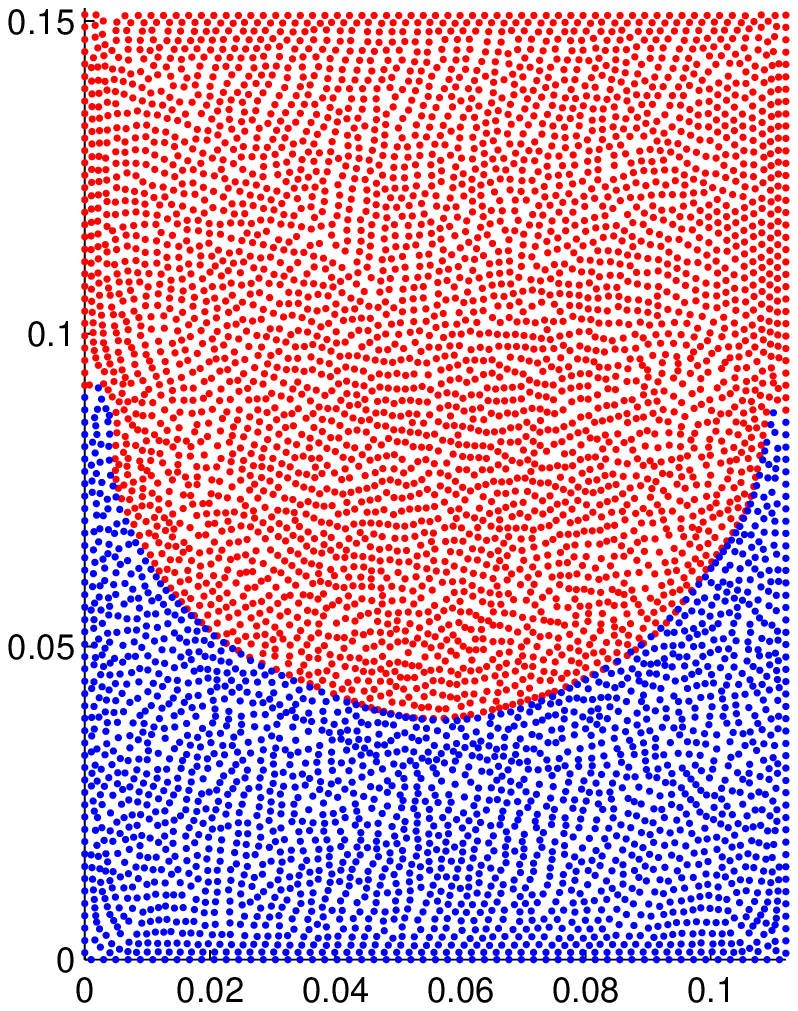} 
 \caption{ Time evolution of the two phases for $\theta_s = 5^{o}$. 
  From left to right, first row: $t=0.5$ s and $1.0$ s, second row: $t=1.5$s and $2.0$ s, 
 third row: $t=2.4$s and $3.0$ s.
  Blue (or dark gray) particles indicate liquid particles and red (or light gray) ones indicate gas particles. }
  \label{pool_theta_5}
\end{figure}

\subsection{Stationary shape of sessile drops} 
\subsubsection{Drops without gravity}
\label{test3} 
This test case is also taken from the paper by Liu et al \cite{LKO}. Consider a two-dimensional domain $[0,0.4]\times[0,0.12]$. The initial 
spacing of the particles is defined via $\Delta x \approx h/3 $ with $h = 0.005$. This gives an initial total number of particles equal to $6536$. 
We consider an initially rectangular ethanol drop of size $[0.15,0.25] \times [0,0.06]$ at the bottom of the domain. If the 
initial particle positions lie inside this rectangle, we define them as liquid particles, otherwise as gas particles, see Fig. \ref{rectangular_drop}(a). 
The fluid parameters are same as in \cite{LKO}. 
The ethanol drop has a density $\rho_l=797.88~ \mbox{kg m}^{-3}$ and a viscosity $\mu_l = 0.0018$ Pa~s. 
The surface tension coefficient is $\sigma = 0.02361 ~\mbox{N m}^{-1}$. 
The gas density is $\rho_g = 1~ \mbox{kg~m}^{-3}$. To dampen the drop oscillations, a large gas viscosity $\mu_g = 0.01$ Ps~s is considered. Gravitational 
and other external forces are set to zero. The initial velocity of all particles is zero. We apply no-slip boundary conditions at all domain boundaries. 
We have studied five different static contact angles  $\theta_s = 30^{o}, 60^{o}, 90^{o}, 120^{o}$ and $150^{o}$. In Fig. \ref{rectangular_drop} the 
stationary drop shapes are shown for all contact angles. These results are comparable with the results presented in Liu et al  \cite{LKO}.   
\begin{figure}
 \includegraphics[keepaspectratio=true, width=.5\textwidth]{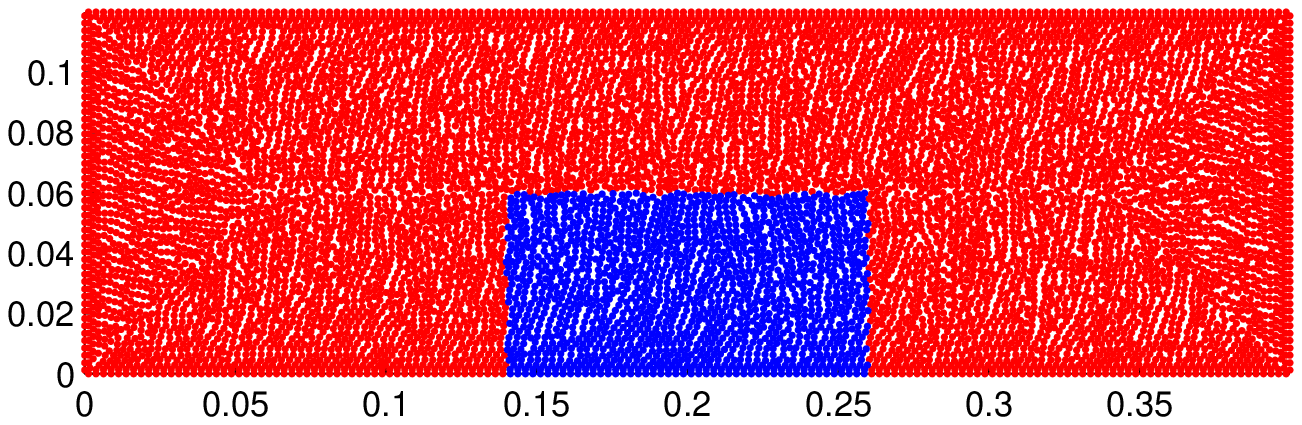}
 \includegraphics[keepaspectratio=true, width=.5\textwidth]{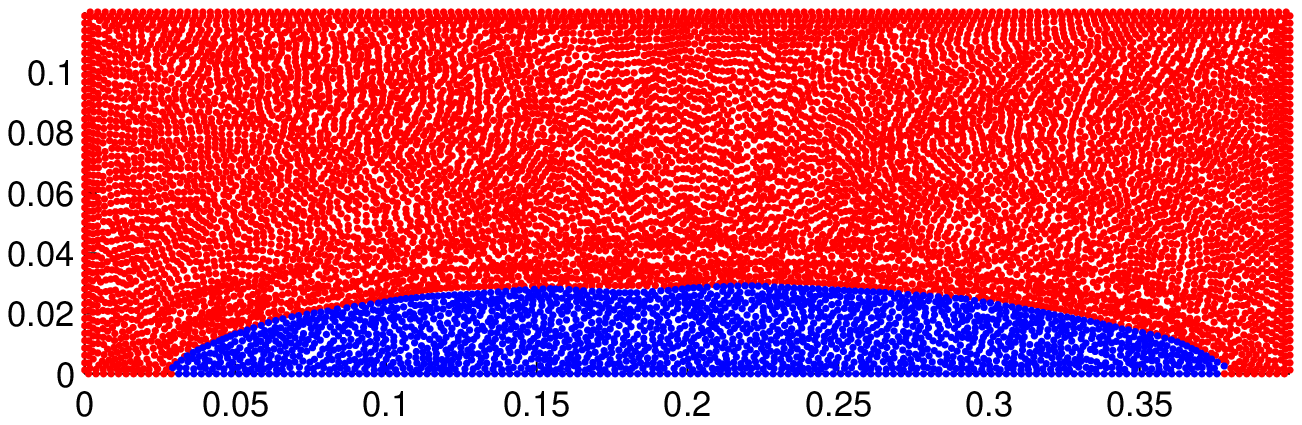}
\\
 \includegraphics[keepaspectratio=true, width=.5\textwidth]{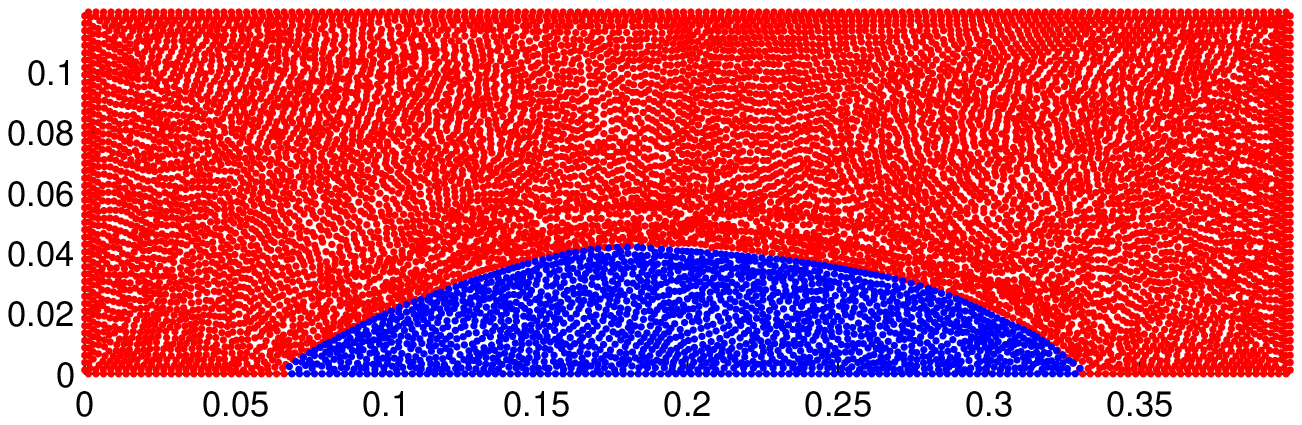}
\includegraphics[keepaspectratio=true, width=.5\textwidth]{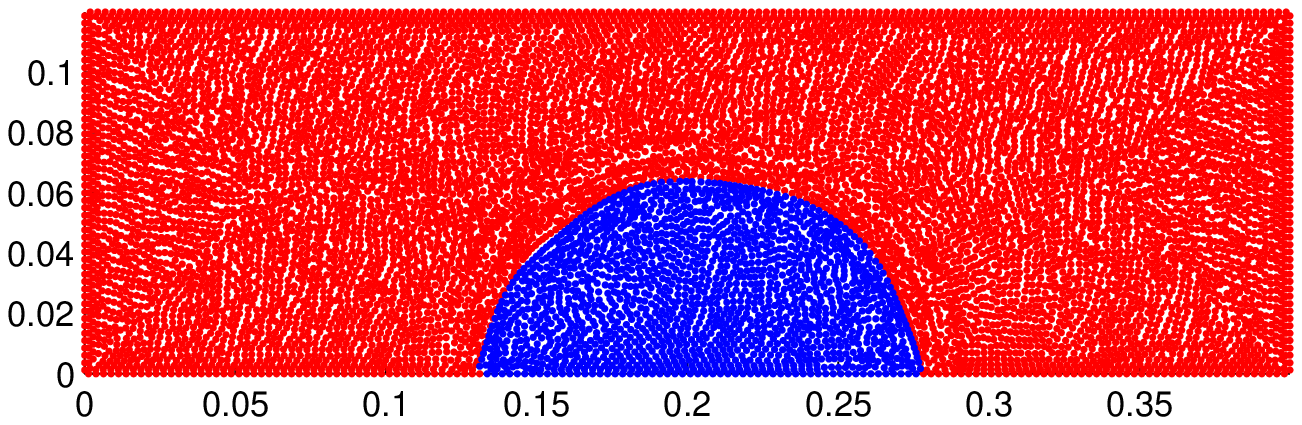}
\\
 \includegraphics[keepaspectratio=true, width=.5\textwidth]{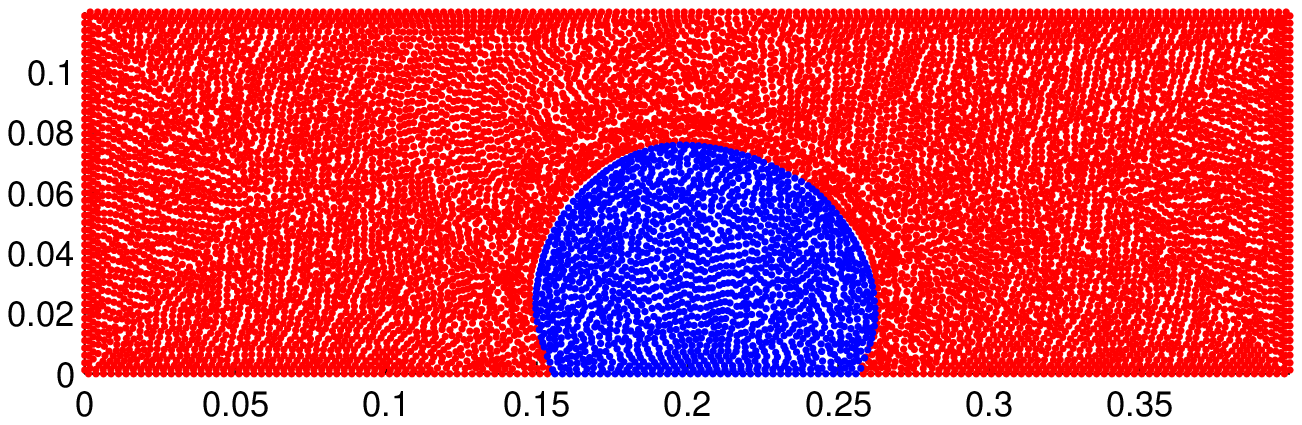}
 \includegraphics[keepaspectratio=true, width=.5\textwidth]{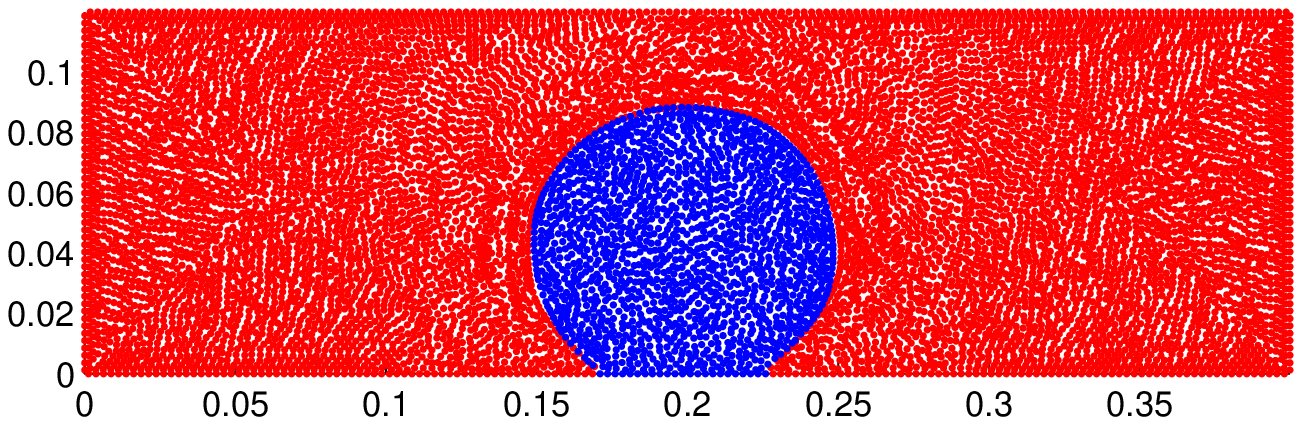} 
 \caption{ Initial shape and stationary shapes of sessile drops.  
  From left to right, first row: initial drop and stationary drop for $\theta_s = 30^{o}$ , second row: $\theta_s = 60^{o}$ and $\theta_s = 90^{o}$, 
 third row: $\theta_s=120{^o}$ and $\theta_s = 150^{o}$. 
 Blue (or dark gray)  indicate liquid particles and red (or light gray) indicate gas particles. }
 \label{rectangular_drop}
\end{figure}

In order to verify the the computational results quantitatively, we consider a circular drop sitting at the bottom wall. We initialize an ethanol drop with radius $R_0$ as shown in Fig. \ref{analytical_drop}. This type of drop was considered in \cite{DL}.

\begin{figure}[h!]
		\begin{tikzpicture}[scale=0.9]
			\draw (-3,0)--(3,0);
			\draw (-2,0) arc (180:0:2cm);
			\draw (0,0)--(1,1.8) node[midway,right=0.1cm]{$R_0$};
 		\end{tikzpicture}
		 \begin{tikzpicture}[scale=0.9]
			\draw (-3,-0)--(3,-0) node[midway,below=0.05cm]{$L$};
			\draw (-2,-1) arc (180:0:2cm);
			\draw (0,0)--(0,1) node[midway,left=0.05cm]{$H$};
			\draw (0,-1)--(1.8,0) node[midway,below=0.05cm]{$R$};
			\draw (-1,0) arc (0:70:0.5cm) node[below=0.01cm] {$\theta_s$};
 		\end{tikzpicture}
	 \caption{ Initial (left) and final (right) shape of a circular drop. }
 \label{analytical_drop}
\end{figure}
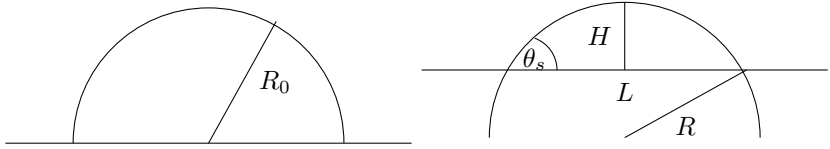

First, we neglect external forces such as gravity.  In a stationary state a circular drop with static contact angle $\theta_s$ will form a circular cap with 
 spreading length $L$, height $H$, and radius of curvature $R$, as shown in Fig. \ref{analytical_drop}. These parameters depend on the drop volume (i.e. area in 2D) and the static contact angle $\theta_s$.  
 The height and spreading length of the drop can be expressed as 
 \begin{equation}
 H = R(1-\mbox{cos}\theta_s), \quad \quad L = 2R\mbox{sin}\theta_s. 
 \label{HandL_of_drop}
 \end{equation}
 From simple geometrical considerations the area of the circular cap in Fig. \ref{analytical_drop} is given by 
 \begin{equation}
 A = R^2\theta_s - R^2 \mbox{sin}\theta_s\mbox{cos}\theta_s. 
 \end{equation}
 Conservation of volume yields $A = \pi^2 R_0^2/2$. So the final radius of curvature is given by 
 \begin{equation}
 R = R_0\sqrt{\frac{\pi}{2(\theta_s - \mbox{sin}\theta_s \mbox{cos}\theta_s )}}. 
 \label{radius_of_drop}
 \end{equation}
 Moreover, in steady state the Laplace law 
 \begin{equation}
 \Delta P = \sigma/R
 \label{laplace_law}
 \end{equation}
 holds, where $\Delta P $ is the pressure difference between the surrounding gas and  the liquid drop.

 We consider a two-dimensional rectangular domain $[0,0.3]\times[0,0.123]$. We initialize liquid particles inside the semicircle of radius $0.06$ and center 
 $(0.15,0)$, and gas particles in the rest of the domain. All fluid parameters are same as in the previous case with an initially rectangular drop except for the viscosities. We have considered larger viscosities 
 such that the fluids reach their equilibrium states faster. Specifically, we have chosen $\mu_g = 0.01$ Pa~s  and $\mu_l = 0.1 $Pa~s. We consider five different 
 static contact angles $\theta_s = 30^{o}, 60^{o}, 90^{o}, 120^{o}$ and $150^{o}$.  In Fig. \ref{compare_analytical_drop} we plot the analytical circles according to (\ref{radius_of_drop}) and the computed liquid distribution in the stationary state. The plots clearly show that the numerical results are very close to the analytical results. 
 
  \begin{figure}
 \includegraphics[keepaspectratio=true, width=.5\textwidth]{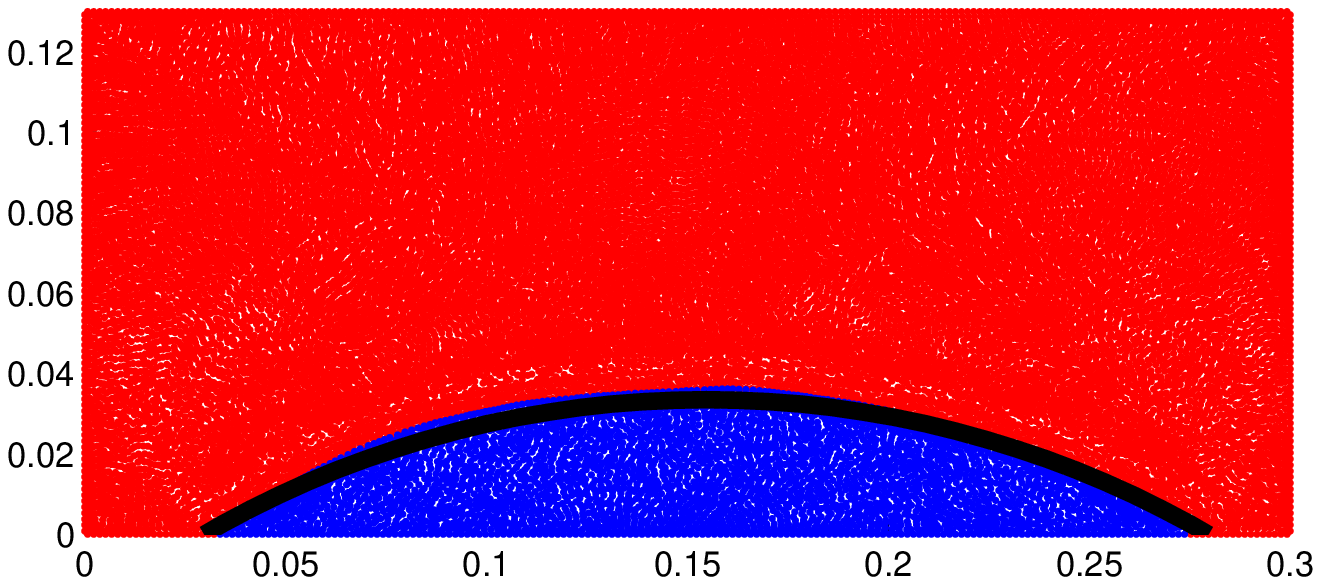}
\includegraphics[keepaspectratio=true, width=.5\textwidth]{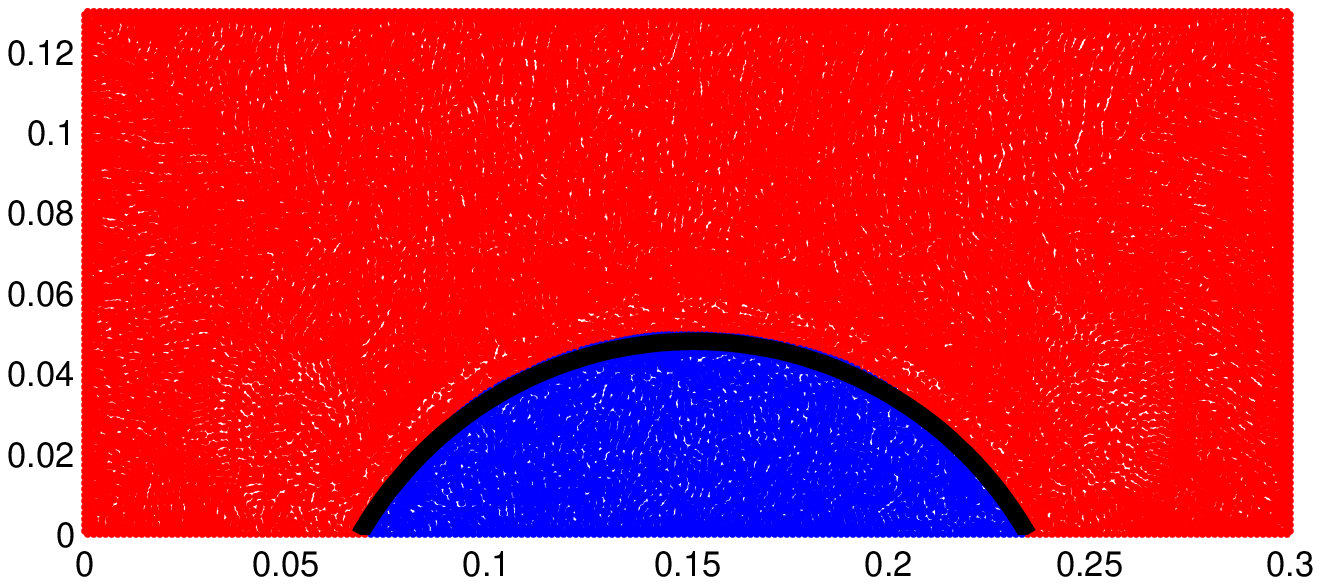}
\\
 \includegraphics[keepaspectratio=true, width=.5\textwidth]{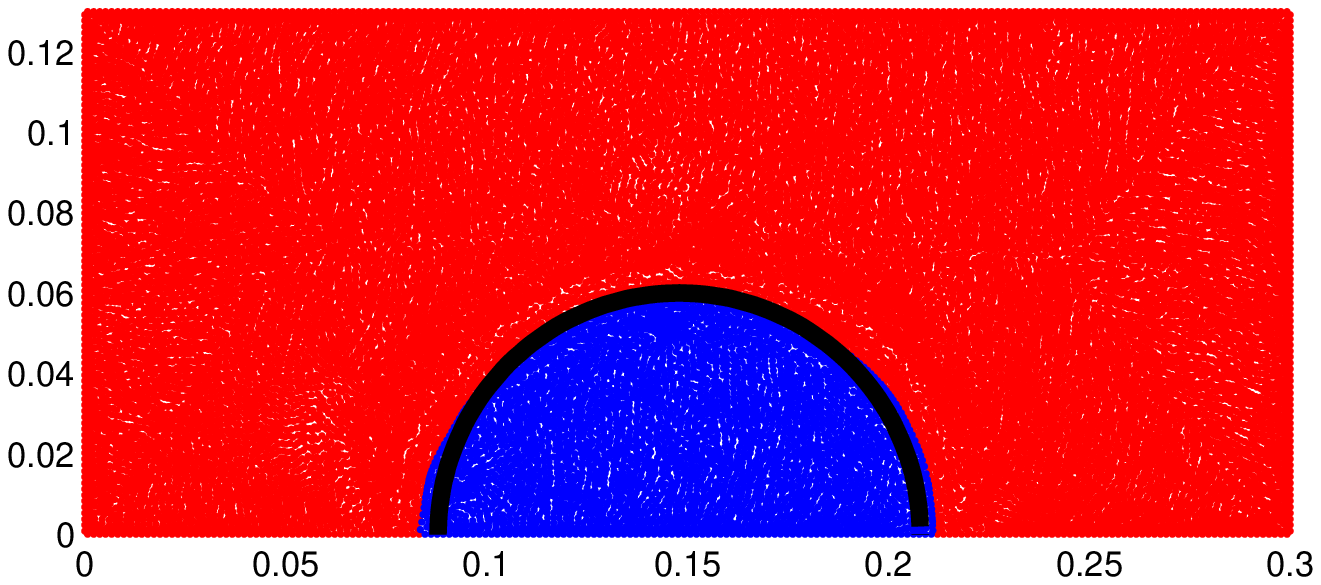}
\includegraphics[keepaspectratio=true, width=.5\textwidth]{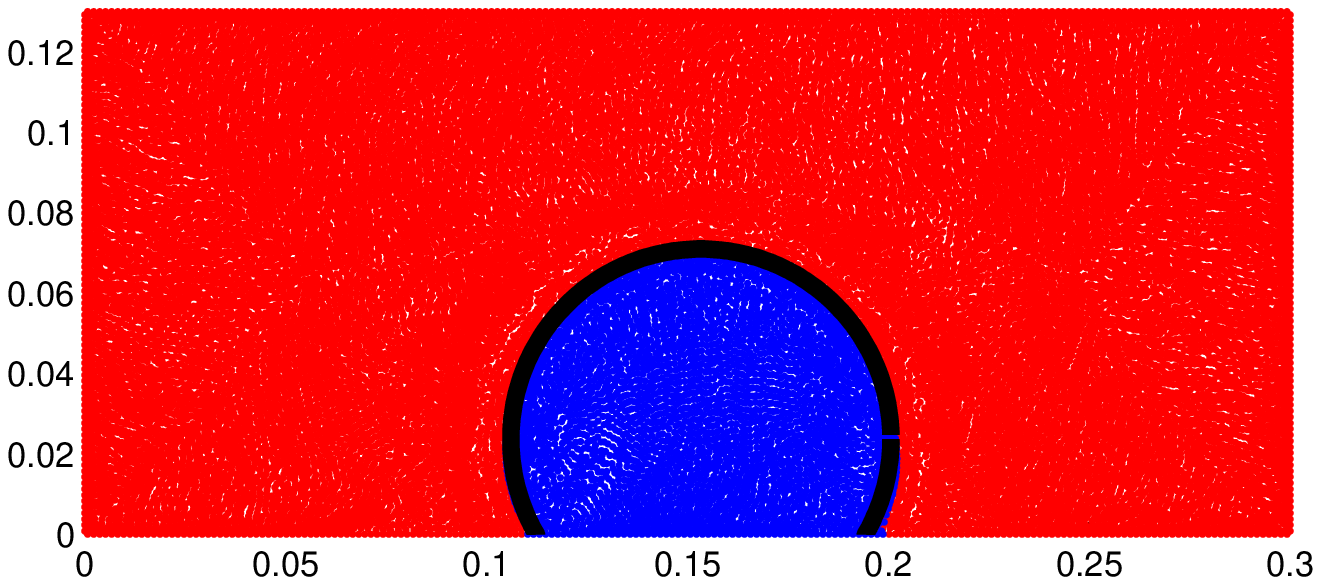}
 \\
 \includegraphics[keepaspectratio=true, width=.5\textwidth]{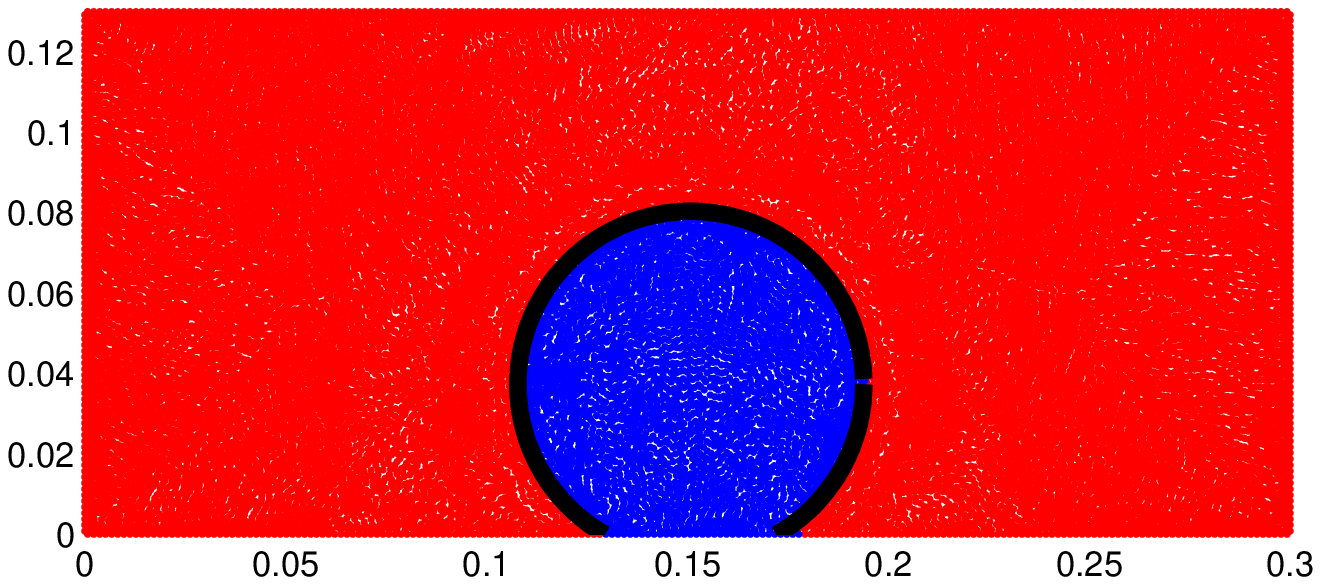} 
 \caption{ Stationary shapes of circular drops. 
 From left to right, first row:  for $\theta_s = 30^{o}$ and $\theta_s = 60^{o}$, second row: $\theta_s = 60^{o}$ and $\theta_s = 90^{o}$, 
 third row: $\theta_s=120{^o}$ and $\theta_s = 150^{o}$. 
 The blue (or dark gray) indicate liquid particles, red (or light gray) indicate gas particles and the solid lines indicate the analytical circular shapes of liquid drops.}
 \label{compare_analytical_drop}
 \end{figure}

Moreover, we plot the analytical spreading length and the height of the circular cap together with their numerical values in Fig. \ref{LandH_of_drop}. 
The numerical value of the spreading length is computed as the difference of the $x$-coordinates of the extreme left and extreme right liquid particles. 
Similarly, the height is computed from the difference of the $y$-coordinates of the extreme top and extreme bottom liquid particles. The results show that for $\theta_s = 30^{o}$ 
it takes longer to converge to the analytical solution, however, for $\theta_s = 150^{o}$ the numerical solution reaches the stationary configuration after 5 seconds. 
Furthermore, we plot the Laplace law (\ref{laplace_law}) in Fig. \ref{pressure_difference}, giving quite satisfactory results. 
The pressure difference is calculated as the mean  pressure of all liquid particles minus the mean pressure of all gas particles. To demonstrate the stability of the scheme, we have computed the kinetic energy inside the liquid drop. In Fig. \ref{kin_energy} we see that after $t=10s$ the kinetic energy remains 
small and stable.    

\begin{figure}
\centering
\includegraphics[keepaspectratio=true, width=.45\textwidth]{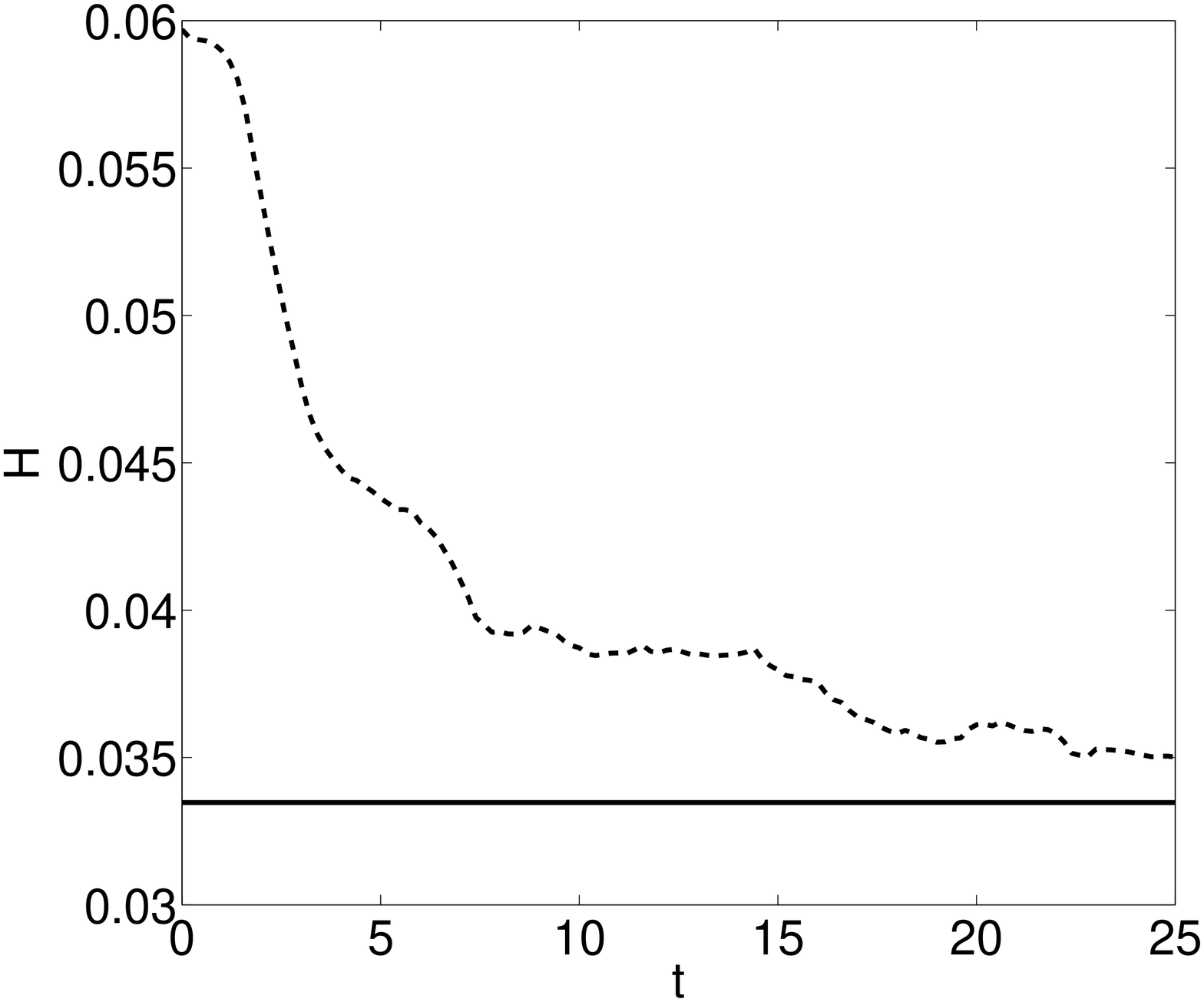}
\includegraphics[keepaspectratio=true, width=.45\textwidth]{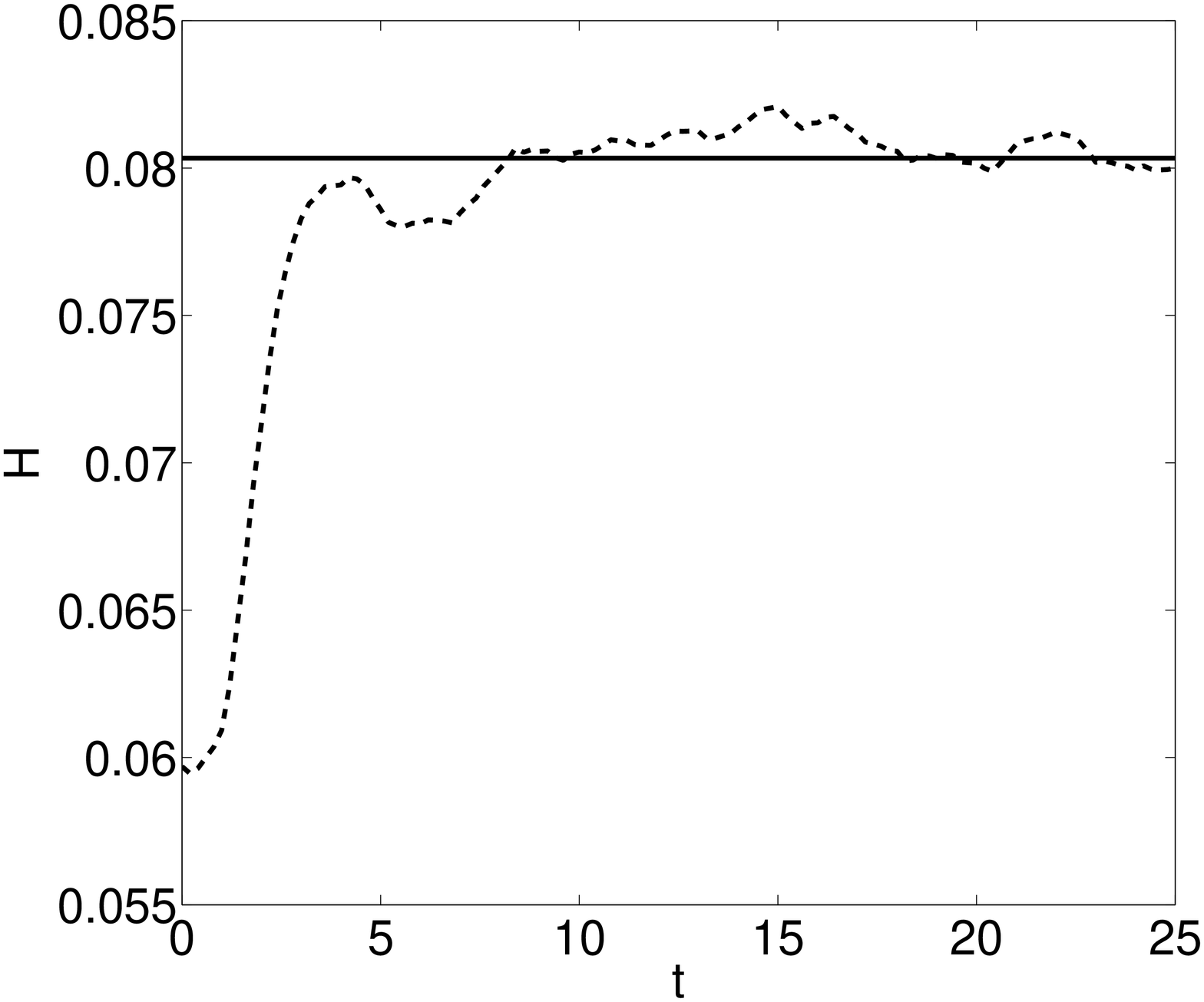}
\\
\includegraphics[keepaspectratio=true, width=.45\textwidth]{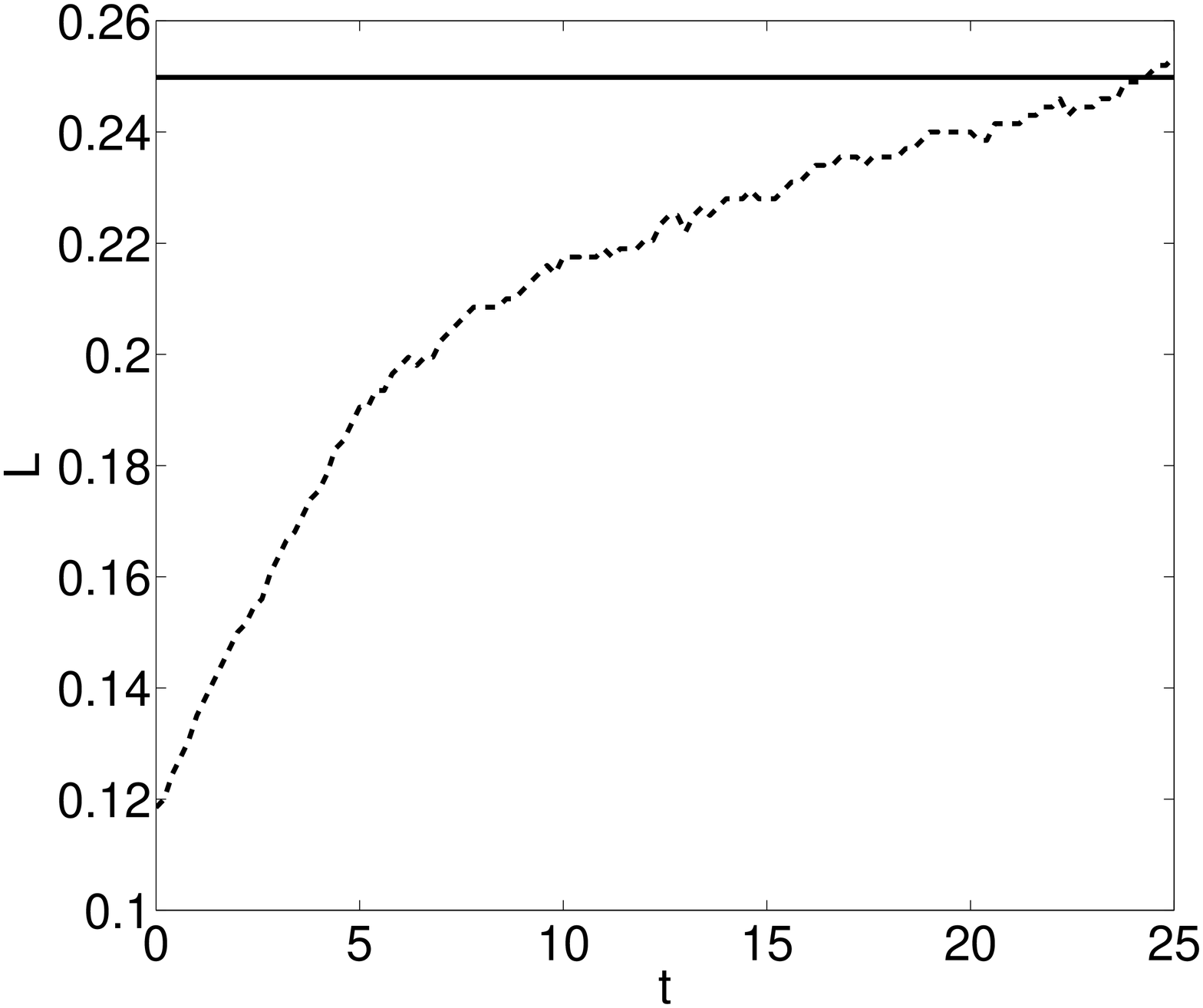} 
\includegraphics[keepaspectratio=true, width=.45\textwidth]{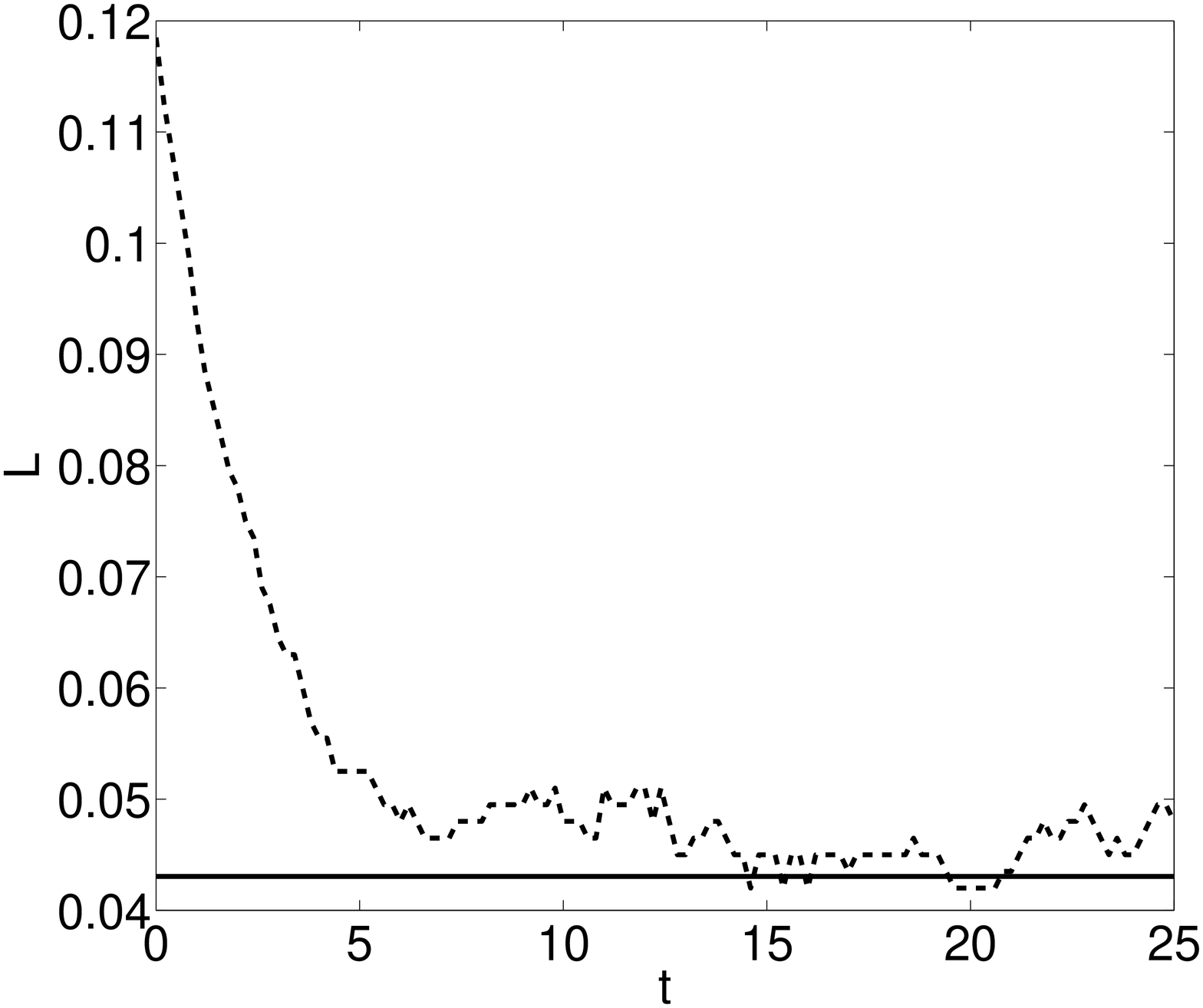}
\caption{ Time evolution of the drop spreading length $L$  and the height $H$. The left figures are for $\theta_s = 30^{o}$, those on  the right  for $\theta_s=150^{o}$. Solid lines  represent the analytical solutions and dash lines represent the numerical solutions.}
\label{LandH_of_drop}
 \end{figure}
 
\begin{figure}
\centering
\includegraphics[keepaspectratio=true, width=.45\textwidth]{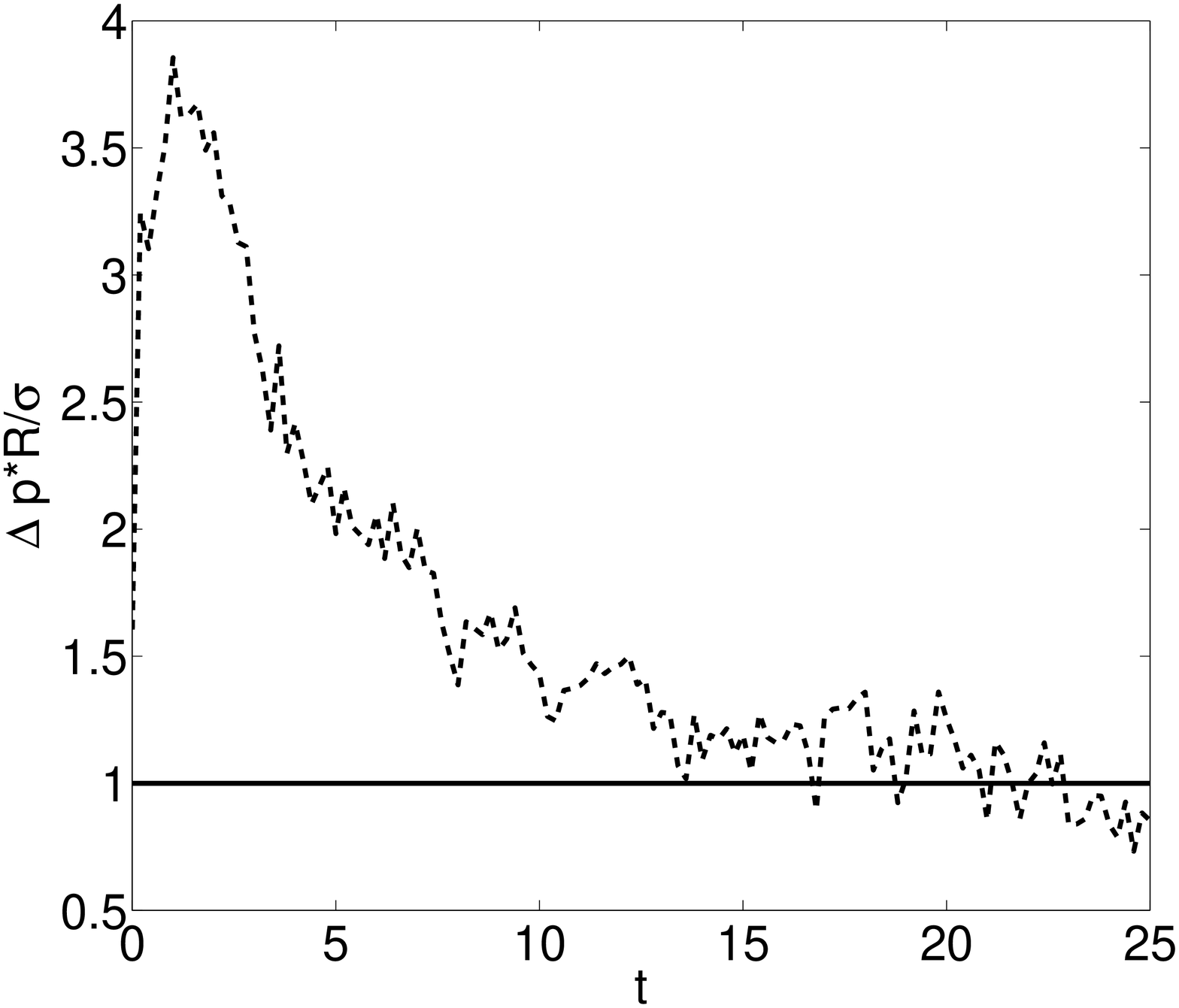}
\includegraphics[keepaspectratio=true, width=.45\textwidth]{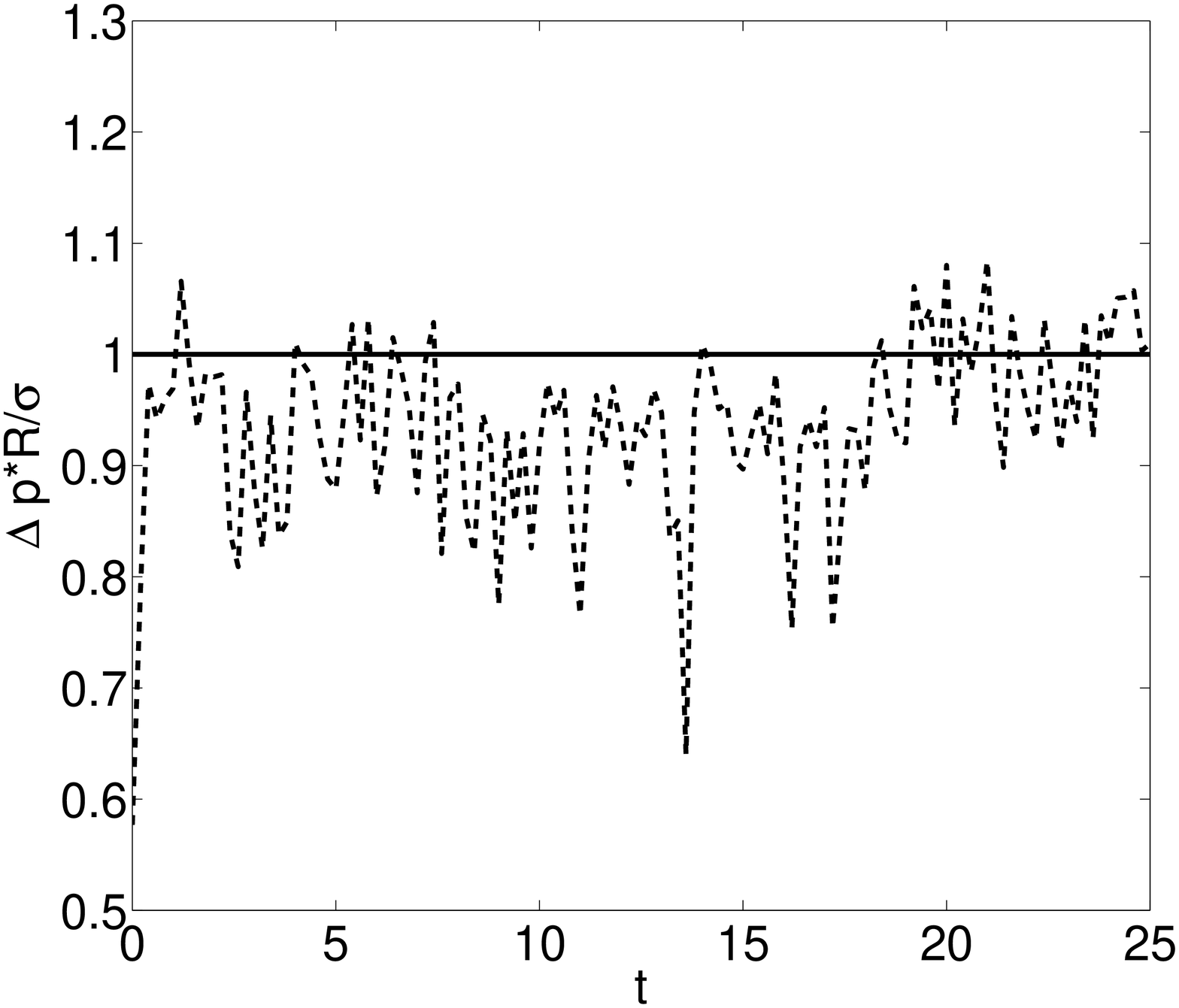}
 \caption{$\Delta p*\frac{R}{\sigma}$ versus time. The left figure is for $\theta_s = 30^{o}$, the right  one for  $\theta_s=150^{o}$. Solid line represent the 
analytical solutions and dash lines represent the numerical solutions.}
 \label{pressure_difference}
  \end{figure}
 
\begin{figure}
  \centering
\includegraphics[keepaspectratio=true, width=.5\textwidth]{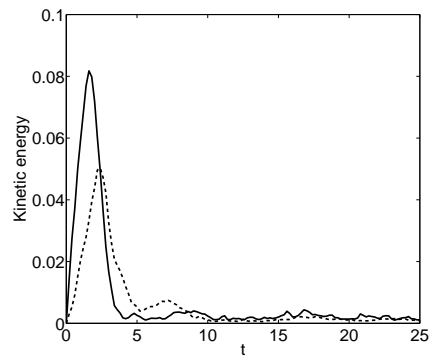}
\caption{Kinetic energy versus time inside the drop. Dash line represents for $\theta_s=30^{o}$ and solid line represents for $\theta_s = 150^{o}$. } 
\label{kin_energy}
\end{figure}

\subsubsection{Gravity effect}
\label{test4}
In this subsection we study the flattening of the drop due to gravity. In the absence of gravity the final shape of drop depends only on the static contact angle $\theta_s$. When gravity is included, the shape also depends on the E\"otv\"os number ($Eo$) defined by 
$Eo = \rho_l g R_0^2/\sigma$,  where $R_0$ is the initial radius, $g$ the gravitational accelaration and $\rho_l$ the liquid density. Sessile drops under gravity have been 
studied in \cite{DL}, where the authors have computed the asymptotic shape of drops as a function of $Eo$. For $Eo << 1$ 
the shape is dominated by surface tension and the drop resembles a circular cap with angle $\theta_s$. In this case the height of the drop can be reexpressed from 
(\ref{HandL_of_drop}) and (\ref{radius_of_drop}) as 
\begin{equation}
H_0 = R_0(1-\cos\theta_s)\sqrt{\frac{\pi}{2(\theta_s - \mbox{sin}\theta_s\mbox{cos}\theta_s)}}.
\label{asymptoticH0}
\end{equation}
For $Eo>>1$   the shape of the drop is dominated by gravity and similar to a pancake whose height is directly proportional to the capillary length \cite{DL}
\begin{equation}
H_{\infty} = 2 \sqrt{\frac{\sigma}{\rho_l g} } \mbox{sin}(\frac{\theta_s}{2}). 
\label{asymptoticHinf}
\end{equation}
 \begin{figure}
 \centering
\includegraphics[keepaspectratio=true, width=.5\textwidth]{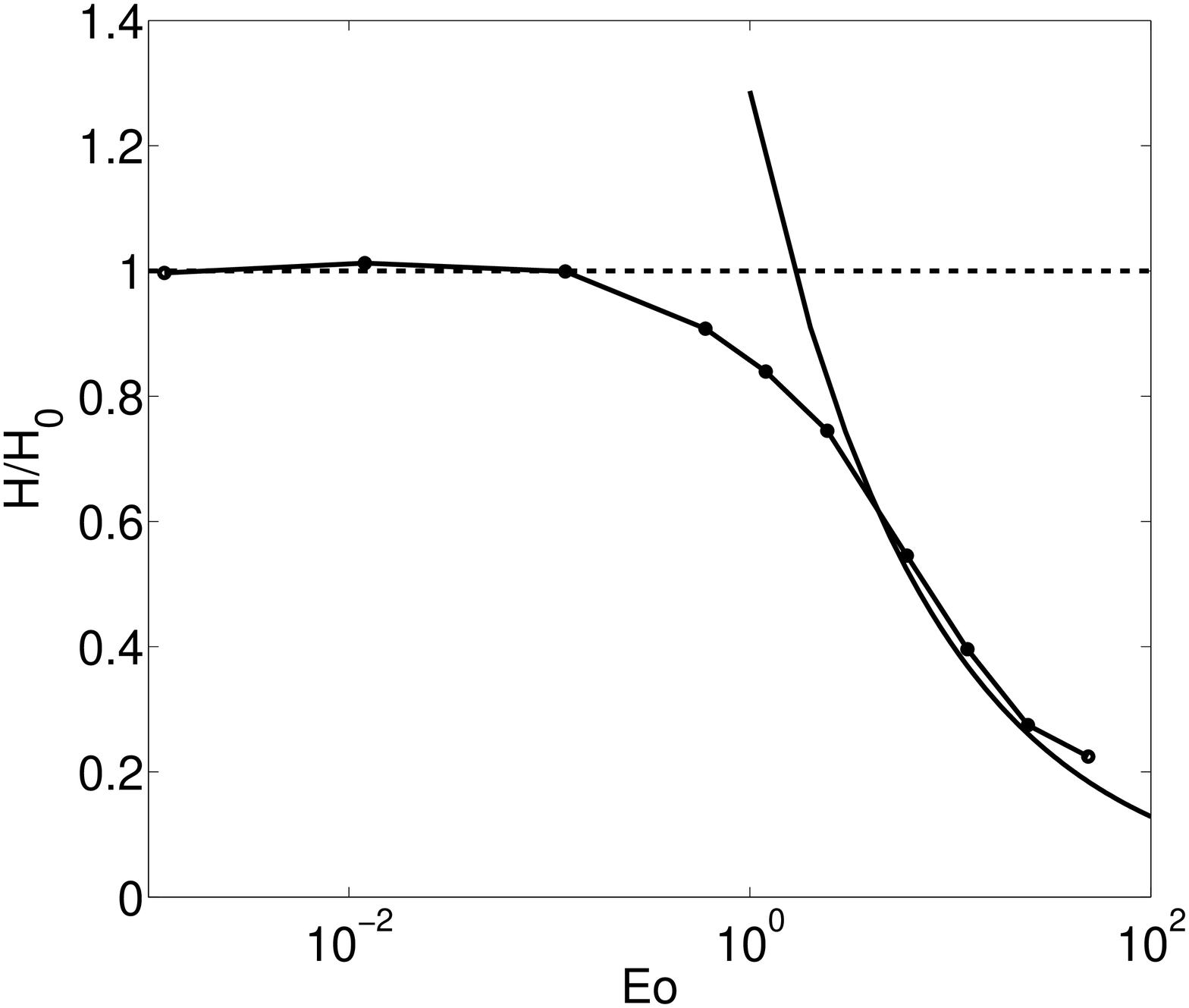}
\caption{ The ratio $H/H_0$ as a function of $Eo$ for $\theta_s = 130^{o}$. Line with '-o-' represents the numerical results,  dash line represents the result of Eq.(\ref{asymptoticH0}), and the solid line  represents the result of  Eq. (\ref{asymptoticHinf}). }
\label{H_vs_Eo}
\end{figure}
 
We again consider an ethanol drop. All initial configurations and parameters are the same as in the previous test case \ref{test3}, except for the contact angle $\theta_s$. Here $\theta_s = 130^{o}$ is chosen in order to compare to the results presented in \cite{DL}.  We change the E\"otv\"os number by changing the gravitational acceleration, keeping the other parameters fixed. The gravity vector points in negative $y$ direction. 
For the cases with $Eo > 12$ we have doubled the width of the 
computational domain which is now given by $[0,0.6] \times[0,0.123]$, since the length wetted by the droplet increases. Then the center of the drop is located at $(0.3,0)$. 
As in the earlier case, we define the height $H$ of the drop as the maximum value of the $y-$coordinates of all liquid particles. As can be seen in Fig. \ref{LandH_of_drop}, the 
height and spreading length of the drop fluctuate around their analytical values. In order to obtain the height $H$ as a function of $Eo$, we take the average value 
with respect to time. The averaging starts at $t = 4s$ and ends at the final time  $t=15s$. 
In Fig. \ref{H_vs_Eo} we show the numerical values of the ratio of $H/H_0$ as a function of E\"otv\"os number 
together with the asymptotic solutions  (\ref{asymptoticH0}) and (\ref{asymptoticHinf}).  The numerical and the asymptotic solutions agree very well. 
Like  in the article by Dupont and Legendre \cite{DL}, we also observe the transition between the circular cap and the pancake shape at around $Eo = 1$. 
In Fig. \ref{drop_with_gravity} we show the final shapes of the drop for varying E\"otv\"os numbers.

\begin{figure}
\centering
\includegraphics[keepaspectratio=true, width=.45\textwidth]{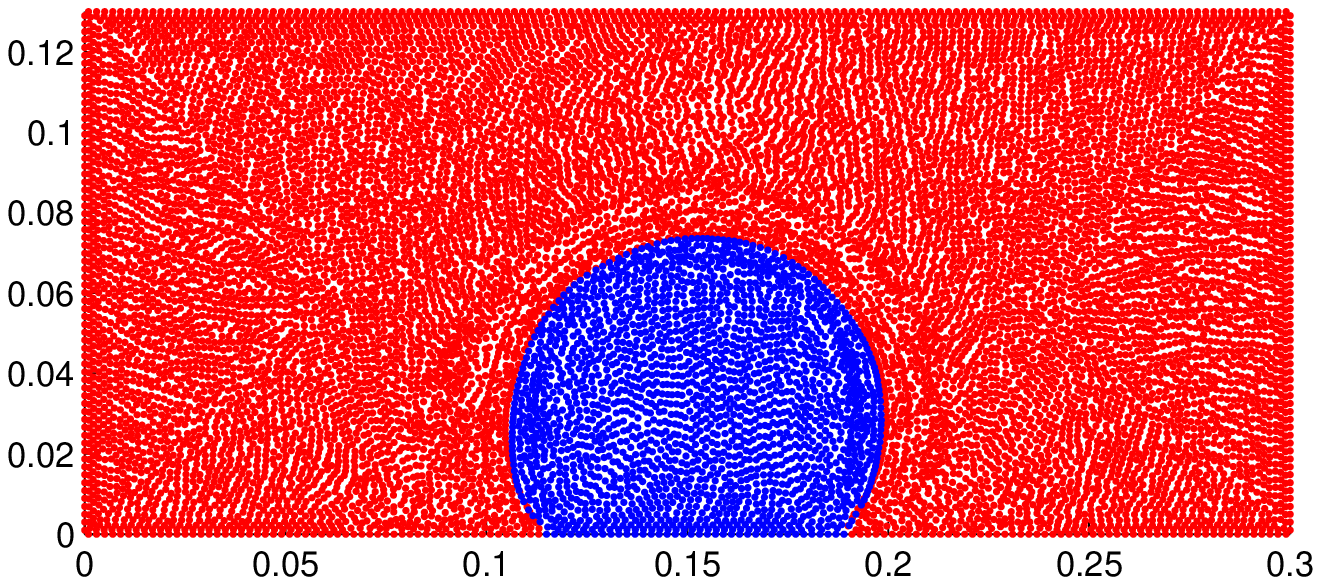}
\includegraphics[keepaspectratio=true, width=.45\textwidth]{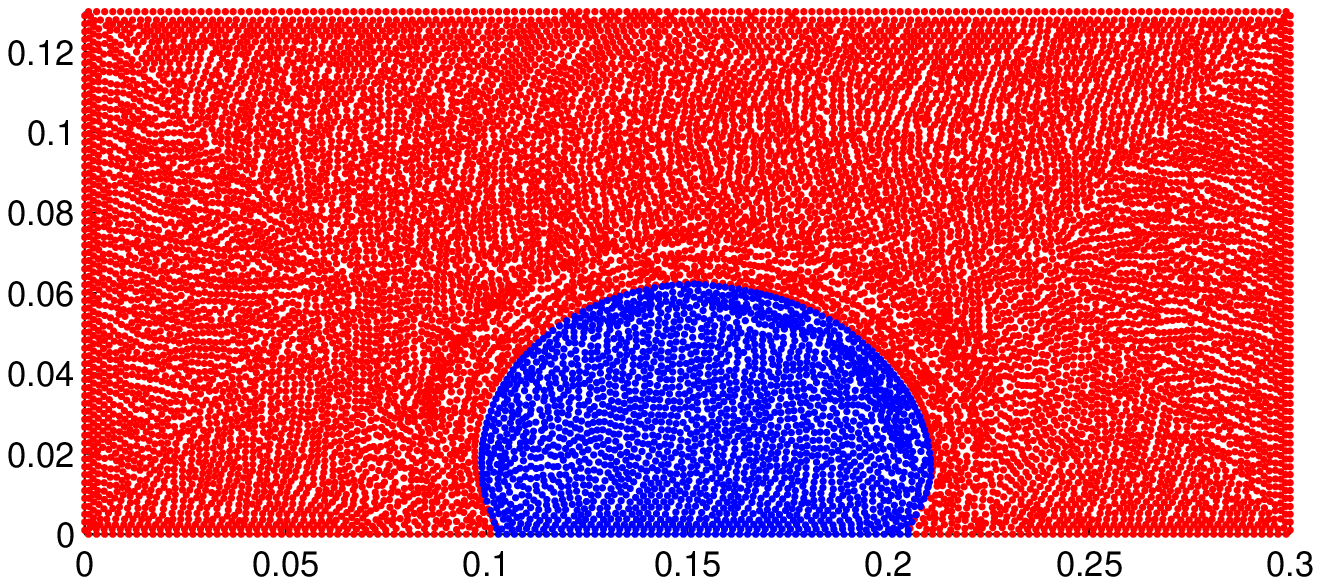}
\\
\includegraphics[keepaspectratio=true, width=.45\textwidth]{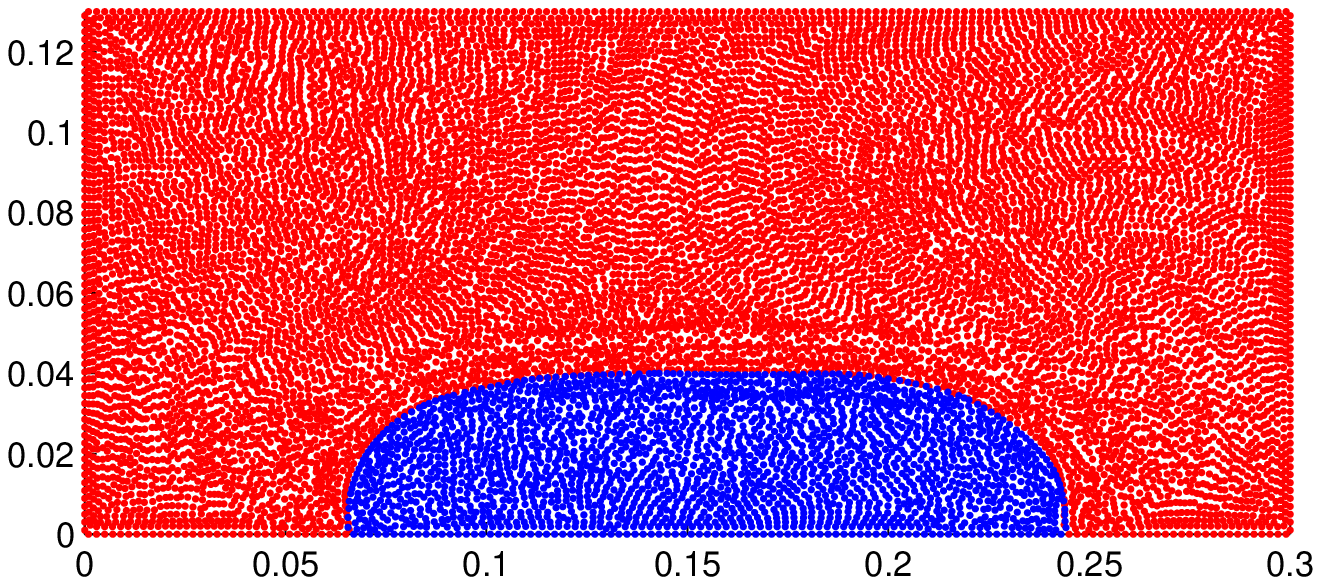} 
\includegraphics[keepaspectratio=true, width=.45\textwidth]{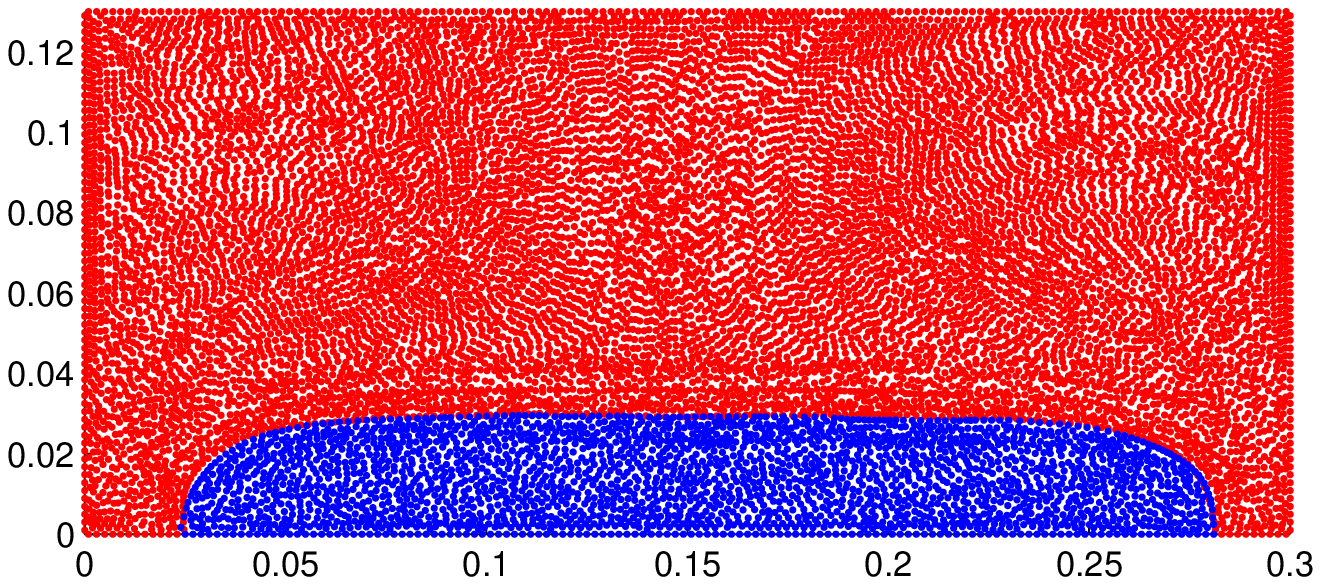}
\caption{ Final shapes of sessile drops for different E\"otv\"os numbers and a static contact angle of $\theta_s = 130^{o}$,  where red (or light gray) indicates gas particles and blue (or dark gray) indicates liquid particles. First row: $Eo = 0.12$ (left) and $Eo = 1.2$ (right). Second row: $Eo = 6.06$ (left) and $Eo = 12.16$ (right). }
\label{drop_with_gravity}
 \end{figure}

\subsection{Convergence study}
\label{convergence_study}
In this section we study the numerical convergence. We consider an initially semi-circular drop sitting at the bottom wall, like in Fig. \ref{analytical_drop}(a). 
The static contact angle is chosen as $\theta_s=150^{o}$. We study three different values of the scale for the interaction radius in the particle method, $h = 0.1, 0.05$ and $0.0025$. The corresponding initial numbers of 
particles are $2186, 8375$ and $32729$, respectively. The other parameters are same as in subsection \ref{test3}. In Fig. \ref{convergence_study} we show 
the values of $R$, $L$ and $\Delta p*R/\sigma$ as a function of time. One can observe convergence of  the numerical solutions to the analytical ones when the number of particles 
increases. For $h=0.05 $ and $0.0025 $ there is not much difference in the height of the drop, but still significant differences in the spreading length and the Laplace pressure are visible between these two cases. 

\begin{figure}
\centering
\includegraphics[keepaspectratio=true, width=.45\textwidth]{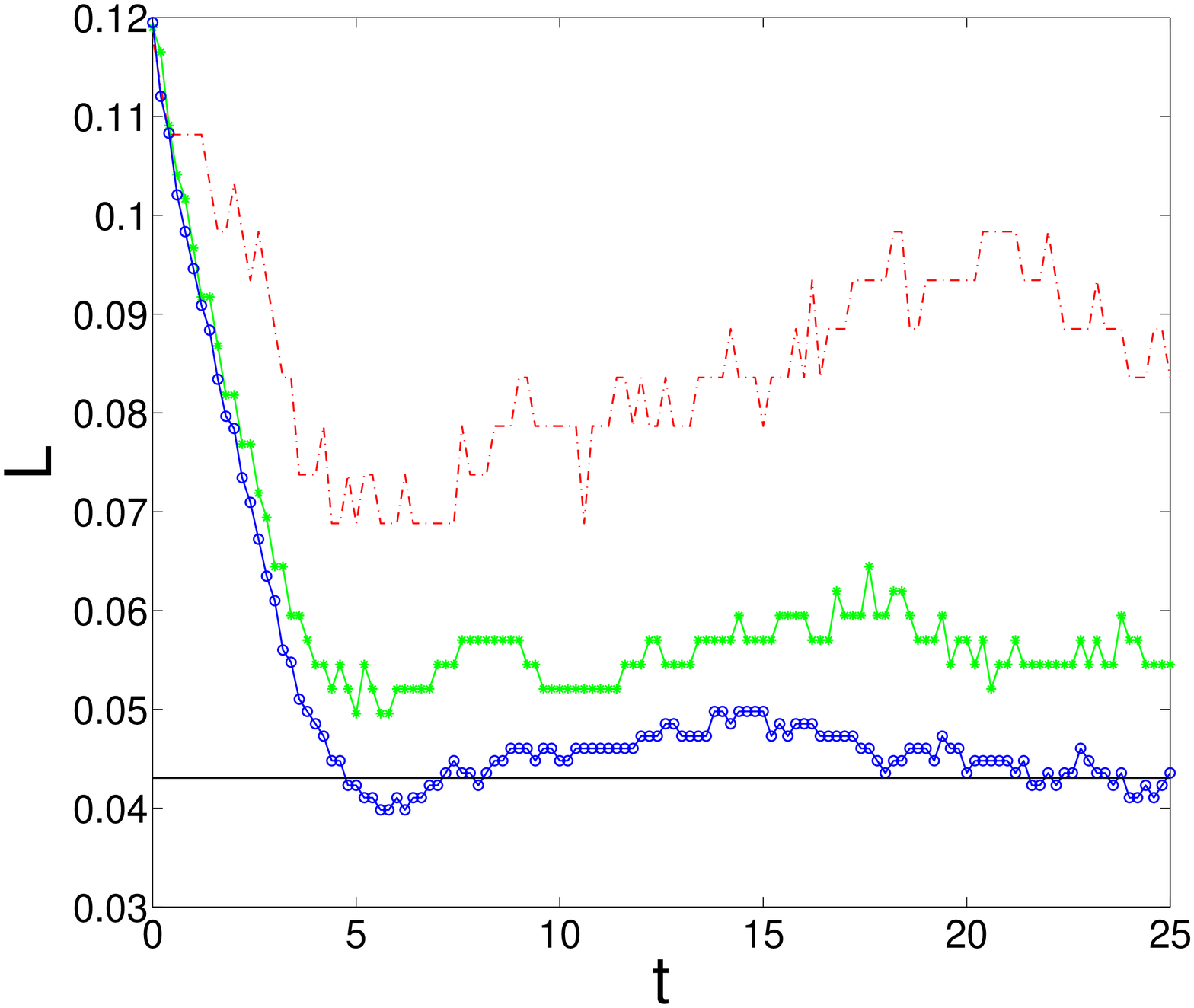}
\hfill
\includegraphics[keepaspectratio=true, width=.45\textwidth]{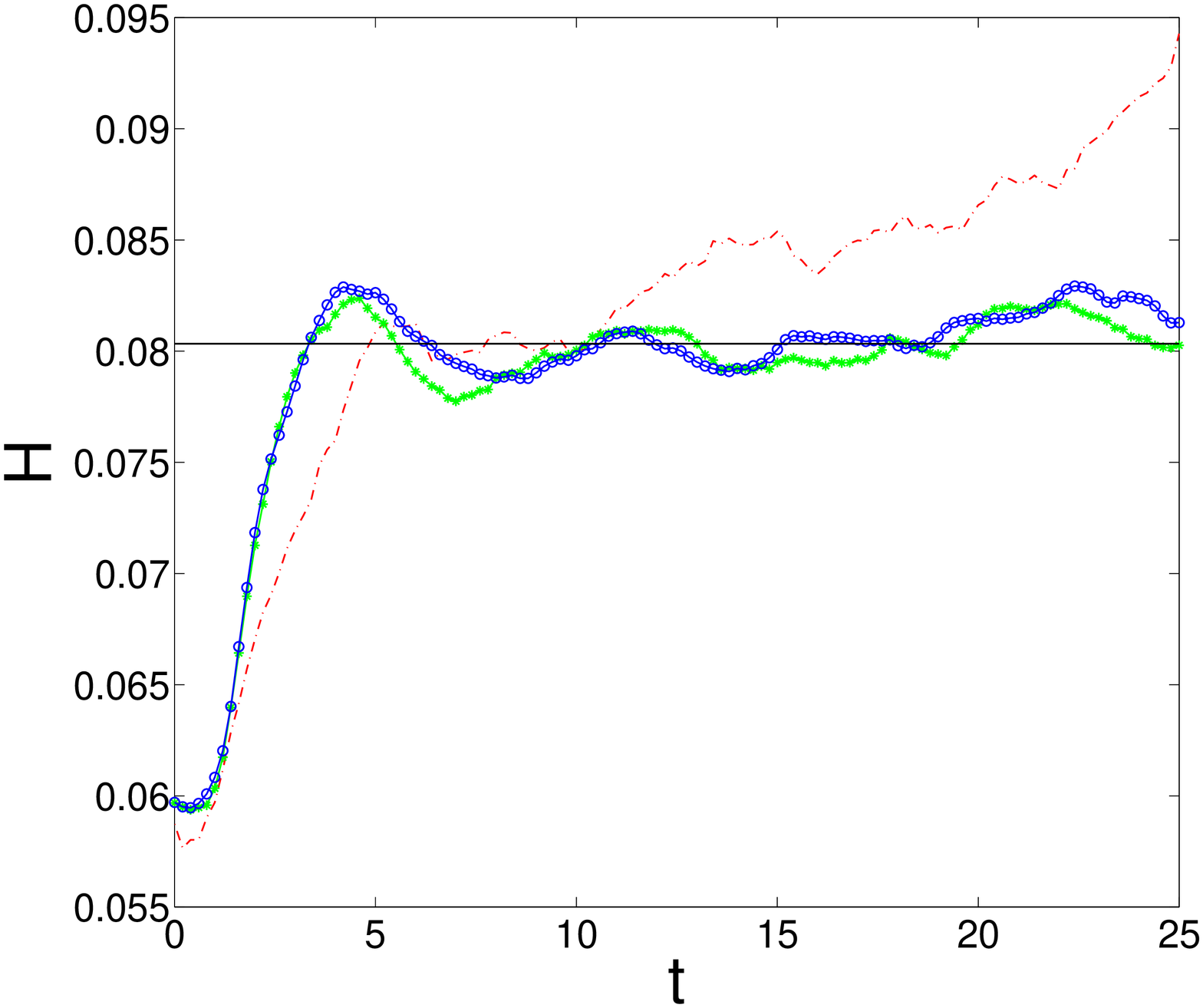}
\\
\centering
\includegraphics[keepaspectratio=true, width=.45\textwidth]{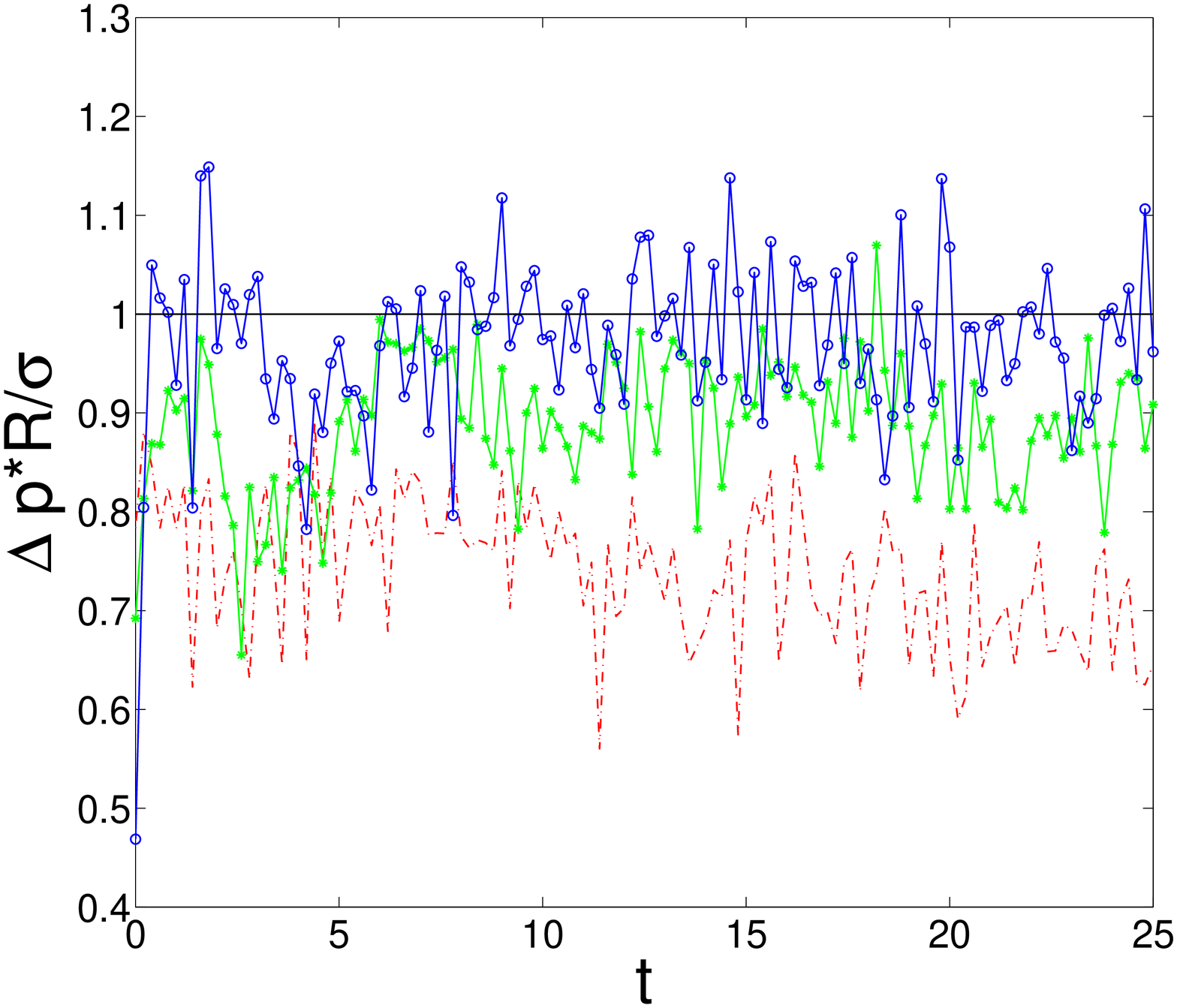} 
 \caption{ Time evolution of the drop spreading length $L$, height $H$  and Laplace pressure for different numbers of particles and $\theta_s = 150^{o}$. Solid lines are 
 analytical values, numerical solutions with red lines $-\cdot-$ are for $h=0.01$, with green lines $-*-$  for $h=0.005$ and with blue lines $-o-$ for $h=0.0025$. }
\label{convergence_study}
 \end{figure}

\section{Concluding Remarks}
\label{conclusion}

We have presented a meshfree Lagrangian particle method to compute two-phase flows driven by wetting forces. The continuous surface force model is used to implement surface tension forces. 
The normal vector and the curvature of the fluid-fluid interface is computed with the help of a color function indicating the two phases based on a least-squares approach. All 
differential operators as well as the solution of the pressure Poisson equation are computed via the least-squares method. At the three-phase contact line between the two fluids and the solid wall the static contact angle is prescribed to model different wetting properties. The numerical results were compared with specific analytical solutions of a diffusion equation with phase-dependent diffusion coefficient. In addition to that, numerical test cases for wetting in a container as well as sessile drops with and without gravity were studied. Apart from a few exceptions where asymptotic analytical results are available, these cases can only be accessed numerically, but have been studied before by other authors. In all of the considered scenarios, a good agreement with the benchmark cases was obtained. Future work will be devoted to simulating dynamic wetting processes for which the dependence of the contact angle on the speed of the three-phase contact line will have to be modeled.


\subsection*{Acknowledgment}  This  work is partially supported 
by the German research foundation (DFG) grant number KL
1105/201.  We would like to thank the DFG for the financial support.



\end{document}